\documentclass[11pt,letterpaper]{article}
\usepackage[T1]{fontenc}
\usepackage{lmodern,amsmath,amssymb,amsthm,mathtools}
\usepackage[margin=1in]{geometry}
\usepackage{microtype,enumitem,fancyhdr,needspace}
\usepackage[hidelinks]{hyperref}
\hypersetup{pdftitle={Observable sets for one-dimensional Schrodinger equations with subquadratic and quadratic power potentials},pdfauthor={Yanchuan Wang},pdfsubject={Revised manuscript},pdfkeywords={Schrodinger equation, observable set, window density, minimal observation time}}
\numberwithin{equation}{section}
\newtheorem{theorem}{Theorem}[section]
\newtheorem{lemma}[theorem]{Lemma}
\newtheorem{proposition}[theorem]{Proposition}

\theoremstyle{definition}
\newtheorem{definition}[theorem]{Definition}

\theoremstyle{remark}
\newtheorem{remark}[theorem]{Remark}
\newcommand{\R}{\mathbb R}

\newcommand{\Z}{\mathbb Z}
\newcommand{\ind}{\mathbf 1}
\newcommand{\loc}{\mathrm{loc}}
\newcommand{\obs}{\mathrm{obs}}
\newcommand{\Tmin}{T_{\min,m}}
\newcommand{\norm}[1]{\lVert#1\rVert}

\DeclareMathOperator{\supp}{supp}
\DeclareMathOperator{\Ree}{Re}
\DeclareMathOperator{\Imm}{Im}
\DeclareMathOperator{\osc}{osc}
\title{Observable sets for one-dimensional Schr\"odinger equations\\
with subquadratic power potentials}
\author{Yanchuan Wang}
\date{}
\begin{document}
\maketitle

\begin{abstract}
In this paper, we study observability inequalities for the one-dimensional
Schrödinger equation associated with the operator
$H_m=-\partial_x^2+|x|^{2m}$, where $0<m\leq1$.
For $0<m<1/2$, observability in some positive time, and equivalently
in every positive time, is characterized by the thickness of the
observation set. This extends the characterization of Chen--Liu--Liu at $m=1/2$. For $1/2<m\leq1$, we characterize observable sets in terms
of the asymptotic lower density on the spatial intervals
$I_{m,a}(x)=[x-a|x|^{2m-1},x+a|x|^{2m-1}]$.
Observability in some positive time is equivalent to positive density
at some scale, whereas observability in every positive time is equivalent
to positive density at every scale. Moreover, for $1/2<m<1$, the minimal
observation time is comparable to the square root of the critical
density scale.
\end{abstract}

\tableofcontents
\medskip

\section{Introduction}
\subsection{Problems}
For $0<m\le1$ we study the one-dimensional Schrödinger equation
\begin{equation}\label{equ-sch}
\begin{cases}
i\partial_tu(t,x)=H_mu(t,x), &(t,x)\in(0,\infty)\times\R,\\
u(0,\cdot)=f\in L^2(\R),
\end{cases}
\qquad H_m=-\partial_x^2+|x|^{2m}.
\end{equation}
The operator is the nonnegative self-adjoint operator associated with
\begin{equation}\label{equ-qf}
q_m[u]=\int_\R\bigl(|u'|^2+|x|^{2m}|u|^2\bigr)\,dx,\qquad
Q_m=\{u\in H^1(\R):|x|^mu\in L^2(\R)\}.
\end{equation}
For $f\in L^2(\R)$ the solution is $u(t)=U_m(t)f=e^{-itH_m}f$. All Hilbert spaces are complex. Inner products are
linear in the first argument. We write $|E|$ for the Lebesgue measure of $E$ and $\norm{\cdot}_2$ for the spatial $L^2(\R)$ norm.

\begin{definition}[Observation]\label{def-obs}
Let $E\subset\R$ be measurable. For $T>0$, let $C_{\obs,m}(E,T)$ be the optimal constant in
\begin{equation}\label{equ-obs}
\norm f_2^2\leq C_{\obs,m}(E,T)\int_0^T\norm{\ind_EU_m(t)f}_2^2\,dt
\qquad(f\in L^2(\R)),
\end{equation}
with the convention that $C_{\obs,m}(E,T)=+\infty$ if no finite constant exists.
The set $E$ is some-time observable if $C_{\obs,m}(E,T)<\infty$ for some $T>0$,
and arbitrary-time observable if this holds for every $T>0$.
Set
\begin{equation}\label{equ-mt}
\Tmin(E)=\inf\{T>0:C_{\obs,m}(E,T)<\infty\},\qquad \inf\varnothing=+\infty.
\end{equation}
\end{definition}
By Prouff~\cite[Lemma~A.3]{Prouff}, if $0<\Tmin(E)<\infty$, then $E$ is
not observable in time $\Tmin(E)$. Since the set of observation times is upward closed,
\begin{equation}\label{equ-zt}
\Tmin(E)=0\quad\Longleftrightarrow\quad E\text{ is arbitrary-time observable}.
\end{equation}

\subsection{Main results}
\begin{definition}[Ordinary thickness]\label{def-t}
A set $E$ is $(L,\gamma)$-thick if $L>0$, $0<\gamma\leq1$, and
\begin{equation}\label{equ-t}
|E\cap[x,x+L]|\geq\gamma L\qquad(x\in\R).
\end{equation}
A measurable $E$ is thick if \eqref{equ-t} holds for some $L>0$ and $0<\gamma\le1$.
\end{definition}

\begin{definition}[Density on the force scale]\label{def-d}
For $1/2<m\leq1$, put $s=2m-1\in(0,1]$ and define, for $a>0$, $x\ne0$,
\begin{align}
I_{m,a}(x)&=[x-a|x|^s,x+a|x|^s],\label{equ-id}\\
d_m(E,a)&=\liminf_{|x|\to\infty}
 \frac{|E\cap I_{m,a}(x)|}{2a|x|^s},\label{equ-d}\\
\kappa_m(E)&=\inf\{a>0:d_m(E,a)>0\},\qquad\inf\varnothing=+\infty.
\label{equ-cs}
\end{align}
\end{definition}
By definition, $d_m(E,a)>0$ if and only if there exist $\gamma_a>0$ and $R_a\geq1$ such that
\begin{equation}\label{equ-dd}
|E\cap I_{m,a}(x)|\geq2\gamma_a a|x|^s\qquad(|x|\geq R_a).
\end{equation}
The constants $\gamma_a$ and $R_a$ may depend on $m$, $a$, and $E$. At $m=1$, the defining interval is explicitly
\[
I_{1,a}(x)=[x-a|x|,x+a|x|],\qquad x\ne0.
\]

\begin{theorem}\label{thm-t}
Let $0<m<1/2$ and let $E\subset\R$ be measurable. The following are equivalent:
\begin{enumerate}[label=\textup{(\roman*)},nosep]
\item $E$ is thick.
\item $E$ is observable in some positive time.
\item $E$ is observable in every positive time.
\end{enumerate}
For fixed $m,T,L,\gamma$, there is a finite $C(m,T,L,\gamma)$ such that
$C_{\obs,m}(E,T)\le C(m,T,L,\gamma)$ for every $(L,\gamma)$-thick set $E$.
\end{theorem}

\begin{theorem}\label{thm-d}
Let $1/2<m\leq1$. Then
\begin{align}
E\text{ is some-time observable}
 &\quad\Longleftrightarrow\quad d_m(E,a)>0\text{ for some }a>0,
\label{equ-st}\\
E\text{ is arbitrary-time observable}
 &\quad\Longleftrightarrow\quad d_m(E,a)>0\text{ for every }a>0.
\label{equ-et}
\end{align}
For $1/2<m<1$, there are $0<c_m\leq C_m<\infty$, depending only on $m$, such that
\begin{equation}\label{equ-tc}
c_m\sqrt{\kappa_m(E)}\leq\Tmin(E)\leq C_m\sqrt{\kappa_m(E)}.
\end{equation}
\end{theorem}
\begin{remark}
At $m=1/2$, the characterization in Theorem~\ref{thm-t} is due to
Chen--Liu--Liu~\cite[Theorem~1.2]{CLL}. At $m=1$, positive density at some
scale is equivalent to weak thickness (Lemma~\ref{lem-wt}), that is,
\[
\liminf_{R\to\infty}\frac{|E\cap[-R,R]|}{2R}>0.
\]
\end{remark}

\subsection{Earlier work and proof strategy}
The relation between thickness and the observation of band-limited functions
has classical antecedents in the norm-equivalence results of
Katsnelson~\cite{Kat} and Logvinenko--Sereda~\cite{LS}.
For the heat equation on $\R^n$, Egidi--Veseli\'{c}~\cite{EV} and
Wang--Wang--Zhang--Zhang~\cite{WWZZ} independently characterized
null-controllability by thickness. These results provide a useful comparison
for the geometric conditions considered here.

For $0<m<1/2$, Theorem~\ref{thm-t} characterizes observability by thickness on
intervals of fixed length. For the one-dimensional Schrödinger equation with a real constant potential,
Huang--Wang--Wang
\cite[Theorem~1.1]{HWW} obtained this criterion for observability in some positive
time. This includes the free Schrödinger equation when the potential is zero.
Su--Sun--Yuan \cite[Theorem~1.2]{SSY} proved observability in every positive
time from thick sets for bounded continuous real potentials. Thickness also enters
the quantitative Logvinenko--Sereda theorem \cite{Kov} and its application to
fractional Schrödinger equations \cite{MPSgeom}. In Euclidean space, other
sufficient conditions involve periodic open sets \cite{Taufer} or a dynamical
condition on the observation region \cite{GK}. For surveys of Schrödinger control
and the geometry of observable sets, see \cite{Laurent,WangICM}.

For observations at isolated times, Wang--Wang--Zhang~\cite{WWZ}
established observability and unique continuation inequalities for the free
Schr\"odinger equation, including a connection between two-time observation
and an uncertainty principle. Huang--Soffer~\cite{HS} developed an abstract
approach to two-time observation using an operator uncertainty principle
and minimal velocity estimates. These time-sampling results differ from
the time-integrated observation considered in \eqref{equ-obs}.
On rectangular two-dimensional tori, Burq--Zworski~\cite{BZ}
proved observability from every measurable spatial set of positive measure,
illustrating the dependence of the geometric criterion on the spatial setting.

The closest comparison for $0<m<1/2$ is the linear potential at $m=1/2$.
For $V(x)=|x|$, Huang--Wang--Wang \cite[Theorem~1.3]{HWWlinear} combined an
Ingham-type estimate with a quantitative compactness argument to bound the
observation cost when $|E^c\cap[x,x+1]|=O(|x|^{-\alpha})$ for some $\alpha>1/2$.
Chen--Liu--Liu \cite[Theorem~1.2]{CLL} obtained the full thickness criterion.
We adapt their resolvent estimate near the turning region and their time-filtering
argument \cite[Sections~3--7]{CLL}. Two features of $|x|^{2m}$ enter the proof.
Its derivative is singular at zero, so we first work with the Lipschitz potential
$\max\{1,|x|^{2m}\}$ and absorb the bounded difference in the resolvent estimate.
At centers $x_n\to\pm\infty$, the translated potential
$|x_n+y|^{2m}-|x_n|^{2m}$ converges locally uniformly to zero as $n\to\infty$, whereas
$|x_n+y|-|x_n|$ converges to $\pm y$. Thus the limiting equation is free,
both in the compactness argument and in the construction that proves necessity.

For $1/2<m\leq1$, we begin with the harmonic oscillator at $m=1$.
Huang--Wang--Wang \cite[Theorem~1.2]{HWW} proved that $E$ is observable in
some positive time if and only if it is weakly thick. The uniform gap between
consecutive eigenvalues allows them to reduce observability in some positive
time to a uniform lower bound on the $L^2(E)$-mass of the normalized
eigenfunctions \cite[Proposition~3.3]{HWW}. They use the parity and asymptotic
behavior of these eigenfunctions to identify this condition with weak
thickness. By Lemma~\ref{lem-wt}, weak thickness is equivalent to
$d_1(E,a)>0$ for some $a>0$.
For classical nonharmonic Fourier estimates under spectral gap assumptions,
including the treatment of finitely many exceptional frequencies, see
Ingham~\cite{Ingham} and Haraux~\cite{Haraux}.

For $1/2<m<1$, the consecutive gaps satisfy
$\lambda_{n+1}-\lambda_n\asymp_m r_n^{m-1}\to0$ as $n\to\infty$, where
$r_n=\lambda_n^{1/(2m)}$. A spectral interval of fixed width can therefore
contain increasingly many eigenvalues. Bounds on individual eigenfunctions
alone do not control cancellation among their linear combinations on $E$.
Throughout $1/2<m\leq1$, we use the spectral gap and pointwise eigenfunction
estimates of Dall'Ara--Martini \cite{DM1,DM2} to control the mass of the
corresponding spectral projection near and beyond the turning points.
In the allowed region, a local Prüfer-angle estimate controls the mass on
each interval by its observed mass and a residual term. Combining these
estimates at a scale $a$ with $d_m(E,a)>0$ gives
\begin{equation}\label{equ-ir}
\norm v_2^2\leq N_a\norm{\ind_Ev}_2^2+
A_ma\norm{(H_m-\lambda)v}_2^2,
\qquad\lambda\geq\Lambda_a,\quad v\in D(H_m).
\end{equation}
Here $N_a$ and $\Lambda_a$ may depend on $m,a,E$, while $A_m$ depends only
on $m$. At $m=1$, a term of order $\sqrt a$ in the projection estimate
persists as $\lambda\to\infty$. We absorb it by restricting to $0<a\leq a_0$.
Miller's time-Fourier argument \cite[Theorem~5.1]{Miller} then yields
observability for $T>\pi\sqrt{A_ma}$.

Conversely, we test observability with a Gaussian packet initially at rest
at $R$. Its displacement over the observation interval is bounded by
$2mT^2R^{2m-1}$ for sufficiently large $R$, with $T\leq1/4$ at $m=1$.
Estimating the time spent by the packet near each point of $E$, we deduce
that observability in time $T$ implies $d_m(E,(2m+1)T^2)>0$ in this range. Together with \eqref{equ-ir}, this proves that
observability in every positive time is equivalent to $d_m(E,a)>0$ for every
$a>0$ throughout $1/2<m\leq1$. For $1/2<m<1$, the same bounds give
$\Tmin(E)\asymp_m\sqrt{\kappa_m(E)}$. Prouff's general sufficient condition
for smooth subquadratic confining potentials \cite[Theorem~1.3]{Prouff}
uses a spatial enlargement of the observation set. Our density condition is
imposed directly on $E$.

\subsection{Declaration of AI Use}
During the preparation of this manuscript, the author used ChatGPT (OpenAI)
to discuss mathematical ideas and proof strategies, refine the presentation
of mathematical arguments, improve the language, and assist with LaTeX
formatting. The author has checked the mathematical derivations, proofs,
and references and takes full responsibility for the content of the final
manuscript.

\section{Observability for \texorpdfstring{$0<m<1/2$}{0 < m < 1/2}}
\subsection{Preliminaries}\label{sec2-pre}
For spatial functions $u(x)$, the form domain $Q_m$ is complete for
$(\norm u_2^2+q_m[u])^{1/2}$: a Cauchy sequence converges in $H^1$ and
in the weighted $L^2$ space, and these limits agree locally.
Smooth cutoffs followed by mollification show that $C_c^\infty(\R)$ is
form-dense. The representation theorem for closed forms \cite{Davies} gives
\begin{equation}\label{equ-dom}
D(H_m)=\{u\in Q_m:-u''+|x|^{2m}u\in L^2(\R)\text{ in distributions}\}.
\end{equation}
For the nontrivial inclusion, test the distributional identity against
$C_c^\infty(\R)$ and extend it to $Q_m$ by form-norm density.
Local boundedness of the potential gives $D(H_m)\subset H^2_{\loc}(\R)$.
Thus $u,u'$ have continuous local representatives.

The embedding $Q_m\hookrightarrow L^2(\R)$ is compact.
On bounded intervals this is Rellich compactness; outside them,
\begin{equation}\label{equ-ft}
\int_{|x|>R}|u|^2\,dx\leq R^{-2m}\norm{|x|^mu}_2^2.
\end{equation}
A diagonal local subsequence and this uniform tail bound give global
compactness. Consequently $H_m$ has compact resolvent and an orthonormal
eigenbasis. If the zero set of an eigenfunction has positive measure, it has an accumulation
point at which both the eigenfunction and its derivative are zero.
Uniqueness for the second-order ODE then gives the zero solution.

\begin{lemma}[A cutoff identity]\label{lem-ci}
Let $H_V$ be the form operator of a nonnegative locally bounded potential $V$,
and let $\eta$ be a real bounded Lipschitz function. For every $u\in D(H_V)$ and $\lambda\in\R$,
\begin{equation}\label{equ-ci}
\int_\R\bigl(|(\eta u)'|^2+(V-\lambda)|\eta u|^2\bigr)\,dx
=\Ree\langle(H_V-\lambda)u,\eta^2u\rangle
 +\int_\R|\eta'|^2|u|^2\,dx.
\end{equation}
\end{lemma}
\begin{proof}
Set $g=(H_V-\lambda)u$. Since $\eta^2u$ belongs to the form domain,
testing against $\eta^2u$ and taking real parts gives
\[
\Ree\langle g,\eta^2u\rangle
=\Ree\int_\R u'\overline{(\eta^2u)'}\,dx
+\int_\R(V-\lambda)|\eta u|^2\,dx.
\]
Since $\eta$ is real,
\[
\Ree\bigl(u'\overline{(\eta^2u)'}\bigr)
=\eta^2|u'|^2+2\eta\eta'\Ree(u'\overline u)
=|(\eta u)'|^2-|\eta'|^2|u|^2
\quad\text{a.e.}
\]
Substitution in the form identity gives \eqref{equ-ci}.
\end{proof}

A bounded real multiplication potential $B(x)$ gives $D(H_V+B)=D(H_V)$ by the form
characterization.

At $m=1/2$ the bounded-force assertion is \cite[Theorem~1.2]{CLL}.
For the rest of this section fix $p=2m\in(0,1)$ and $V_p(x)=|x|^p$.

For $h\in L^1(\R)$, we use the Fourier transform
\[
\widehat h(\tau)=(2\pi)^{-1/2}\int_{\R}e^{-it\tau}h(t)\,dt.
\]
With $D_t=-i\partial_t$, every $h\in\mathcal S(\R)$ satisfies
\[
\widehat{D_th}(\tau)=\tau\widehat h(\tau),
\qquad
\widehat{i\partial_th}(\tau)=-\tau\widehat h(\tau).
\]

For $\psi\in L^\infty(\R)$ and $h\in L^2(\R)$, define
$\psi(D_t)h=\mathcal F_t^{-1}(\psi\widehat h)$ by Plancherel's theorem.
If $\psi\in C_c^\infty(\R)$, the convolution theorem gives
\[
\begin{aligned}
(\psi(D_t)h)(t)
&=\frac{1}{\sqrt{2\pi}}\int_{\R}e^{it\tau}\psi(\tau)\widehat h(\tau)\,d\tau
=(K_\psi*h)(t),\\
K_\psi(t)
&=\frac{1}{2\pi}\int_{\R}e^{it\tau}\psi(\tau)\,d\tau,
\qquad K_\psi\in\mathcal S(\R).
\end{aligned}
\]

Throughout this section, each $\chi\in C_c^\infty(0,T)$ is extended by zero
to $\R$. If $U$ is defined on $(0,T)$, we use the same convention for $\chi U$.

We use the following three facts from Chen--Liu--Liu \cite{CLL}.

\newcounter{boundedfact}
\refstepcounter{boundedfact}\label{fact-fi}
\noindent\textbf{Fact \theboundedfact. Observation on a fixed interval}
\cite[Lemma~3.4, (3.12)]{CLL}.
Fix $L>0$, $0<\gamma\le1$ and $M\ge0$. There are $k_*>0$ and $C>0$,
depending only on $L,\gamma,M$, with the following property. Let $I$ be
an interval of length $L$, let $S\subset I$ be measurable with
$|S|\ge\gamma L$, and let $b\in L^\infty(I;\R)$ satisfy $\norm b_\infty\le M$.
For every $k\ge k_*$ and $u\in H^2(I)$,
\[
\norm u_{L^2(I)}^2\le C\norm u_{L^2(S)}^2+
\frac C{k^2}\norm{-u''+(b-k^2)u}_{L^2(I)}^2.
\]
This follows by translating $I$ to $(0,L)$ and squaring the cited estimate.

\refstepcounter{boundedfact}\label{fact-fc}
\noindent\textbf{Fact \theboundedfact. Compactness after time filtering}
\cite[Lemma~2.8, (2.14)--(2.15)]{CLL}.
Let $T>0$ and let $\rho$ be a positive measurable weight with
$\rho,\rho^{-1}\in L^\infty_{\loc}(\R)$. Suppose that $U_n\rightharpoonup U$
in $L^2((0,T)\times\R,\rho(y)\,dt\,dy)$ as $n\to\infty$ and
\[
i\partial_t U_n=-\partial_y^2 U_n+V_n(y)U_n
\quad\text{in }\mathcal D'((0,T)\times\R),
\]
where $\sup_n\norm{V_n}_{L^\infty(J)}<\infty$ on every bounded interval $J$.
Fix $\chi\in C_c^\infty(0,T)$ and $\psi\in C_c^\infty(\R)$. Along a subsequence,
\[
Z_n:=\psi(D_t)(\chi U_n)\longrightarrow\psi(D_t)(\chi U)
\quad\text{strongly in }L^2_{\loc}(\R_t\times\R_y)
\quad\text{as }n\to\infty.
\]
For each bounded interval $J$, the filtered sequence also satisfies
\[
\lim_{S\to\infty}\sup_n\int_{|t|>S}\int_J|Z_n(t,y)|^2\,dy\,dt=0.
\]

\refstepcounter{boundedfact}\label{fact-fu}
\noindent\textbf{Fact \theboundedfact. Filter removal and unique continuation.}
\begin{enumerate}[label=\textup{(\roman*)},nosep]
\item \emph{Removal of the time multiplier} \cite[Lemma~2.9]{CLL}.
Suppose $g\in L^2(\R)$ has compact support and $\psi$ is bounded and
continuous, with $\psi(\tau)\ne0$ throughout a nonempty open interval.
Then $\psi(D_t)g=0$ implies $g=0$.
\item \emph{Unique continuation from a measurable spatial set}
\cite[Proposition~5.3]{CLL}.
Let $I$ be a connected open time interval, let $V\in L^\infty_{\loc}(\R;\R)$
be independent of time, and let $U\in L^2_{\loc}(I\times\R)$ satisfy
$i\partial_tU=-\partial_y^2U+V(y)U$ in distributions.
If $U=0$ almost everywhere on $I\times S$ for a measurable
$S\subset\R$ of positive measure, then $U=0$ almost everywhere on $I\times\R$.
\end{enumerate}
For unique continuation from open sets for Schrödinger operators with
partially Gevrey coefficients, see \cite{FLL}.

\subsection{A resolvent estimate retaining the turning region}
We adapt the estimate of \cite[Proposition~4.2]{CLL} to the power $p$.
To avoid differentiating $V_p$ at zero, put
\begin{equation}\label{equ-wp}
W_p(x)=\begin{cases}1,&|x|\le1,\\ |x|^p,&|x|>1.\end{cases}
\end{equation}
Since $p|x|^{p-1}\le1$ for $|x|>1$, $W_p$ is $1$-Lipschitz.
Also $0\leq W_p-V_p\leq1$, so
$D(H_{W_p})=D(H_m)$ and $H_{W_p}=H_m+(W_p-V_p)$.

\begin{proposition}\label{prop-tr}
Fix $p\in(0,1)$, $L>0$, and $0<\gamma\leq1$. There exist $C,R_*>0$
such that whenever $R\ge R_*$, every $(L,\gamma)$-thick set $E$ satisfies,
for all $u\in D(H_m)$ and $\lambda\in\R$,
\begin{equation}\label{equ-tre}
\norm u_2^2\leq C\norm{\ind_Eu}_2^2+
\frac CR\norm{(H_m-\lambda)u}_2^2+
6\int_{|W_p-\lambda|\leq2R}|u|^2\,dx.
\end{equation}
\end{proposition}
\begin{proof}
Put $g(x)=[(H_{W_p}-\lambda)u](x)$. Then
\begin{equation}\label{equ-m3}
\begin{aligned}
\norm u_2^2
&=\underbrace{\int_{\{W_p-\lambda\leq-2R\}}|u|^2\,dx}_{\Pi_1}
 +\underbrace{\int_{\{|W_p-\lambda|<2R\}}|u|^2\,dx}_{\Pi_2}
 +\underbrace{\int_{\{W_p-\lambda\geq2R\}}|u|^2\,dx}_{\Pi_3}.
\end{aligned}
\end{equation}

\Needspace{4\baselineskip}\medskip\noindent\textbf{\boldmath Estimate of $\Pi_1$.}\par\nopagebreak[4]
Set
\[
I_j=[jL,(j+1)L),\qquad j\in\mathbb Z,
\]
so $\overline I_j=[jL,(j+1)L]$. The half-open
intervals are pairwise disjoint and cover $\mathbb R$; their closures
intersect only at endpoints. Set
\[
\mathcal J_R=
\left\{j\in\mathbb Z:
\overline I_j\cap
\{x\in\mathbb R:W_p(x)-\lambda\leq-2R\}
\neq\varnothing\right\}.
\]
For each $j\in\mathcal J_R$, define
\[
k_j^2=-\min_{x\in\overline I_j}(W_p(x)-\lambda),\qquad
b_j(x)=W_p(x)-\lambda+k_j^2,\quad x\in\overline I_j.
\]
By the definition of $\mathcal J_R$, we have
$k_j^2\geq2R$, and $b_j\geq0$;
Since $W_p$ is 1-Lipschitz,
\[
0\leq b_j(x)\leq\operatorname{osc}_{\overline I_j}W_p\leq L.
\]
In particular, $\norm{b_j}_{L^\infty(I_j)}\leq L$, and on $I_j$
we have the local identity
\[
g=-u''+(-k_j^2+b_j)u\qquad\text{a.e. on }I_j.
\]
Because $E$ is $(L,\gamma)$-thick, $|E\cap I_j|\geq\gamma L$.
Apply Fact~\ref{fact-fi} on $I_j$ with $S=E\cap I_j$ and $M=L$.
For $R\geq k_*^2/2$, $k_j\ge\sqrt{2R}\ge k_*$. Hence
\[
\begin{aligned}
\int_{I_j}|u|^2\,dx
&\leq C\int_{E\cap I_j}|u|^2\,dx
+\frac{C}{k_j^2}\int_{I_j}|g|^2\,dx\\
&\leq C\int_{E\cap I_j}|u|^2\,dx
+\frac{C}{2R}\int_{I_j}|g|^2\,dx.
\end{aligned}
\]
Summing over $j\in\mathcal J_R$ and using the disjointness of the
intervals gives
\begin{equation}\label{equ-b1}
\begin{aligned}
\Pi_1
&\leq C\sum_{j\in\mathcal J_R}\int_{E\cap I_j}|u|^2\,dx
+\frac{C}{2R}\sum_{j\in\mathcal J_R}\int_{I_j}|g|^2\,dx\\
&\leq C\norm{\ind_Eu}_2^2+\frac{C}{2R}\norm g_2^2.
\end{aligned}
\end{equation}

\Needspace{4\baselineskip}\medskip\noindent\textbf{\boldmath Estimate of $\Pi_3$.}\par\nopagebreak[4]
Define the cutoff by
\[
\eta(x)=
\begin{cases}
0,&W_p(x)-\lambda\leq R,\\[1mm]
\dfrac{W_p(x)-\lambda-R}{R},&R<W_p(x)-\lambda<2R,\\[2mm]
1,&W_p(x)-\lambda\geq2R.
\end{cases}
\]
Since $W_p$ is 1-Lipschitz,
$\eta$ is Lipschitz and
$|\eta'|\leq R^{-1}$ almost everywhere; moreover,
$\eta'=0$ almost everywhere outside
$\{R<W_p-\lambda<2R\}$.
On $\{\eta\neq0\}$ we have $W_p-\lambda>R$. By \eqref{equ-ci},
\[
\begin{aligned}
R\norm{\eta u}_2^2
&\leq\int_\R\bigl(|(\eta u)'|^2+(W_p-\lambda)|\eta u|^2\bigr)\,dx\\
&=\Ree\langle g,\eta^2u\rangle+\int_\R|\eta'|^2|u|^2\,dx\\
&\leq\norm g_2\norm{\eta u}_2
+R^{-2}\int_{R<W_p-\lambda<2R}|u|^2\,dx,
\end{aligned}
\]
where $\norm{\eta^2u}_2\leq\norm{\eta u}_2$ since $0\leq\eta\leq1$.
By Young's inequality,
\[
\norm g_2\norm{\eta u}_2
\leq\frac R2\norm{\eta u}_2^2+
\frac1{2R}\norm g_2^2.
\]
Since $\eta=1$ on $\{W_p-\lambda\geq2R\}$, it follows that
\begin{equation}\label{equ-b3}
\begin{aligned}
\Pi_3&\leq\norm{\eta u}_2^2
\leq R^{-2}\norm g_2^2+
2R^{-3}\int_{R<W_p-\lambda<2R}|u|^2\,dx\\
&\leq R^{-2}\norm g_2^2+2R^{-3}\Pi_2.
\end{aligned}
\end{equation}

\Needspace{4\baselineskip}\medskip\noindent\textbf{Combination of the estimates.}\par\nopagebreak[4]
Substituting \eqref{equ-b1} and \eqref{equ-b3} into
\eqref{equ-m3} gives
\[
\norm u_2^2\leq C\norm{\ind_Eu}_2^2+
\left(\frac{C}{2R}+\frac1{R^2}\right)\norm g_2^2
+\left(1+\frac2{R^3}\right)\Pi_2.
\]
For $R\geq1$, the turning-region coefficient satisfies
$1+2R^{-3}\leq3$, while
$C/(2R)+R^{-2}\leq C_0/R$ with
$C_0=C/2+1$, independent of $R$. Since
$\Pi_2\leq\int_{|W_p-\lambda|\leq2R}|u|^2\,dx$, we obtain
\[
\norm u_2^2\leq C\norm{\ind_Eu}_2^2+
\frac{C_0}{R}\norm g_2^2+
3\int_{|W_p-\lambda|\leq2R}|u|^2\,dx.
\]
Finally, since $0\leq W_p-V_p\leq1$,
\[
\begin{aligned}
\norm g_2^2
&=\norm{(H_{W_p}-\lambda)u}_2^2\\
&=\norm{(H_m-\lambda)u+(W_p-V_p)u}_2^2\\
&\leq2\norm{(H_m-\lambda)u}_2^2+2\norm u_2^2.
\end{aligned}
\]
Substitution yields
\[
\begin{aligned}
\norm u_2^2
&\leq C\norm{\ind_Eu}_2^2+
\frac{2C_0}{R}\norm{(H_m-\lambda)u}_2^2\\
&\quad+\frac{2C_0}{R}\norm u_2^2+
3\int_{|W_p-\lambda|\leq2R}|u|^2\,dx.
\end{aligned}
\]
Set $R_*=\max\{k_*^2/2,1,4C_0\}$. Then $2C_0/R\leq1/2$ for $R\geq R_*$. Absorbing this term into the left-hand side and
multiplying by at most two gives
\[
\norm u_2^2\leq2C\norm{\ind_Eu}_2^2+
\frac{4C_0}{R}\norm{(H_m-\lambda)u}_2^2+
6\int_{|W_p-\lambda|\leq2R}|u|^2\,dx,
\]
which proves~\eqref{equ-tre} after enlarging the constant $C$.
\end{proof}

\subsection{Uniform sufficiency of thickness}\label{sec2-suf}
Fix $p,T,L,\gamma$. The failure of uniform observation gives $(L,\gamma)$-thick sets $E_n$
and $f_n\in L^2(\R)$ such that
\begin{equation}\label{equ-cd}
\norm{f_n}_2=1,\qquad
O_n=\int_0^T\norm{\ind_{E_n}U_m(t)f_n}_2^2\,dt\longrightarrow0
\quad\text{as }n\to\infty.
\end{equation}
By unitarity,
$\norm{\ind_{E_n}U_m(\cdot)h}_{L^2((0,T)\times\R)}\leq\sqrt T\,\norm h_2$
for every $h\in L^2(\R)$, uniformly in $n$.
Since $D(H_m)$ is dense in $L^2(\R)$, we may approximate each $f_n$ there
with $L^2$ error less than $1/(2n)$ and normalize the approximants.
We thus take $f_n\in D(H_m)$ satisfying \eqref{equ-cd}.

\Needspace{4\baselineskip}\medskip\noindent\textbf{Step 1. A lower bound for the mass in the turning region.}\par\nopagebreak[4]
Put $u_n(t,x)=(U_m(t)f_n)(x)$ on $\R^2$.
Choose $I_0=(t_-,t_+)\Subset(0,T)$ and a real time cutoff $\chi(t)$ in
$C_c^\infty(0,T)$ such that $0\leq\chi\leq1$,
$\chi>0$ on $I_0$ and $\chi=0$ outside $I_0$.
So it has support $[t_-,t_+]$.
Set $v_n(t,x)=\chi(t)u_n(t,x)$. Then
\[
i\partial_tv_n=H_mv_n+i\chi'u_n,\qquad
(H_m+\tau)\widehat v_n(\tau)=-i\,\widehat{\chi'u_n}(\tau).
\]
For $R\geq R_*$ apply \eqref{equ-tre} to $\widehat v_n(\tau,\cdot)$
with $\lambda=-\tau$. The transformed equation gives
\[
\begin{aligned}
\norm{\widehat v_n(\tau)}_2^2
&\leq C\norm{\ind_{E_n}\widehat v_n(\tau)}_2^2
+\frac CR\norm{\widehat{\chi'u_n}(\tau)}_2^2\\
&\quad+6\int_{|\tau+W_p(x)|\leq 2R}|\widehat v_n(\tau,x)|^2\,dx .
\end{aligned}
\]
Integrate in $\tau$. By Plancherel in time, $\norm{u_n(t)}_2=1$, and
the commutation of spatial multiplication with the time transform,
\begin{align*}
\int_\R\norm{\widehat v_n(\tau)}_2^2\,d\tau
 &=\int_0^T|\chi(t)|^2\norm{u_n(t)}_2^2\,dt=\norm\chi_{L^2_t}^2,\\
\int_\R\norm{\ind_{E_n}\widehat v_n(\tau)}_2^2\,d\tau
 &=\int_0^T|\chi(t)|^2\norm{\ind_{E_n}u_n(t)}_2^2\,dt\leq O_n,\\
\int_\R\norm{\widehat{\chi'u_n}(\tau)}_2^2\,d\tau
 &=\int_0^T|\chi'(t)|^2\norm{u_n(t)}_2^2\,dt=\norm{\chi'}_{L^2_t}^2 .
\end{align*}
Consequently
\begin{equation}\label{equ-bm}
\norm\chi_{L^2_t}^2\leq CO_n+\frac{C}R\norm{\chi'}_{L^2_t}^2+6M_n,\qquad
M_n=\iint_{|\tau+W_p(x)|\leq2R}|\widehat v_n(\tau,x)|^2\,d\tau\,dx.
\end{equation}
Fix $R\geq R_*$ so the residual term is at most $\norm\chi_{L^2_t}^2/4$.
For all large $n$, $CO_n\leq\norm\chi_{L^2_t}^2/4$, hence
$M_n\geq b_0:=\norm\chi_{L^2_t}^2/12>0$.

\Needspace{4\baselineskip}\medskip\noindent\textbf{Step 2. Selection of the centers $k_n$ and normalization.}\par\nopagebreak[4]
Use the weight $\rho(y)=(1+|y|)^{-4}$ and set
$\mathcal H_\rho=L^2((0,T)\times\R,\rho(y)\,dt\,dy)$.
Let $Q_k=(k-1/2,k+1/2)$.
Choose a frequency cutoff $\psi(\tau)$ in $C_c^\infty(\R)$ with
$0\leq\psi\leq1$ and $\psi=1$ on
$[-2R-2,2R+2]$. Define
\begin{align*}
B_{n,k}&=\norm{\psi(D_t)(\chi e^{itV_p(k)}u_n)}_{L^2(\R\times Q_k)}^2,\\
A_{n,k}&=\int_0^T\int_\R\rho(x-k)|u_n(t,x)|^2\,dx\,dt,\\
D_{n,k}&=\int_0^T\int_\R\rho(x-k)\ind_{E_n}(x)|u_n(t,x)|^2\,dx\,dt.
\end{align*}
On $Q_k$, $|W_p(x)-V_p(k)|\leq|x-k|+1\leq3/2$.
Modulation and Plancherel give
\[
B_{n,k}=\int_{Q_k}\int_\R
|\psi(\tau+V_p(k))|^2|\widehat v_n(\tau,x)|^2\,d\tau\,dx.
\]
The multiplier is one if $|\tau+W_p(x)|\leq2R$, so
$\sum_kB_{n,k}\geq M_n$.
Since
\[
C_\rho:=\sup_{x\in\R}\sum_{k\in\Z}\rho(x-k)<\infty,
\]
Tonelli's theorem and $\|u_n(t)\|_2=1$, $0\leq\chi\leq1$ give
\[
\begin{aligned}
\sum_{k\in\Z}A_{n,k}
&=\int_0^T\int_\R
  \left(\sum_{k\in\Z}\rho(x-k)\right)|u_n(t,x)|^2\,dx\,dt\\
&\leq C_\rho\int_0^T\|u_n(t)\|_2^2\,dt
=C_\rho T.
\end{aligned}
\]
Moreover, for $x\in Q_k$ we have $|x-k|\leq1/2$, and hence
\[
\rho(x-k)=(1+|x-k|)^{-4}\geq(3/2)^{-4}.
\]
Using $0\leq\psi\leq1$ and Plancherel in time, we obtain
\[
\begin{aligned}
B_{n,k}
&=\int_{Q_k}\int_\R
 |\psi(\tau+V_p(k))|^2|\widehat v_n(\tau,x)|^2\,d\tau\,dx\\
&\leq\int_{Q_k}\int_\R|\widehat v_n(\tau,x)|^2\,d\tau\,dx\\
&=\int_{Q_k}\int_0^T|\chi(t)u_n(t,x)|^2\,dt\,dx\\
&\leq(3/2)^4\int_{Q_k}\int_0^T
 \rho(x-k)|u_n(t,x)|^2\,dt\,dx\\
&\leq(3/2)^4 A_{n,k}.
\end{aligned}
\]
Tonelli's theorem also gives
\[
\begin{aligned}
\sum_{k\in\Z}D_{n,k}
&=\int_0^T\int_\R
 \left(\sum_{k\in\Z}\rho(x-k)\right)
 \ind_{E_n}(x)|u_n(t,x)|^2\,dx\,dt\\
&\leq C_\rho O_n\longrightarrow0\qquad\text{as }n\to\infty.
\end{aligned}
\]
The selection lemma \cite[Lemma~6.1]{CLL}, with
$M_0=C_\rho T$, $M_1=(3/2)^4$, $b_0=\norm\chi_{L^2_t}^2/12$,
supplies $\alpha=b_0/(8M_0)>0$ and indices $k_n$ with $A_{n,k_n}>0$ such that
$B_{n,k_n}\geq\alpha A_{n,k_n}$ and
$D_{n,k_n}/A_{n,k_n}\leq\epsilon_n\to0$ as $n\to\infty$.
Put $c_n=\sqrt{A_{n,k_n}}>0$ and
\begin{equation}\label{equ-norm}
U_n(t,y)=c_n^{-1}e^{itV_p(k_n)}u_n(t,k_n+y),\qquad E_n'=E_n-k_n.
\end{equation}
Changing variables gives
\begin{equation}\label{equ-no}
\norm{U_n}_{\mathcal H_\rho}^2=1,\qquad
\int_0^T\int_\R\rho(y)\ind_{E_n'}(y)|U_n(t,y)|^2\,dy\,dt
=\frac{D_{n,k_n}}{A_{n,k_n}}\leq\epsilon_n\longrightarrow0
\quad\text{as }n\to\infty,
\end{equation}
and
\begin{equation}\label{equ-nm}
\norm{\psi(D_t)(\chi U_n)}_{L^2(\R\times Q_0)}^2
=\frac{B_{n,k_n}}{A_{n,k_n}}\geq\alpha>0.
\end{equation}
After translation and normalization, the weighted observed mass tends to zero,
while the filtered mass on $Q_0$ is bounded below by $\alpha$.

\Needspace{4\baselineskip}\medskip\noindent\textbf{Step 3. The limiting equation and a nonzero filtered limit.}\par\nopagebreak[4]
Direct differentiation gives
\[
i\partial_tU_n=-\partial_y^2U_n+V_n(y)U_n,\qquad
V_n(y)=|k_n+y|^p-|k_n|^p.
\]
The inequality $\bigl||a|^p-|b|^p\bigr|\leq|a-b|^p$ gives
$|V_n(y)|\leq|y|^p$.
After extraction the integer centers are constant or tend to one of
the infinities, and $V_n\to V_\infty$ uniformly on compact intervals as $n\to\infty$, where
\[
V_\infty(y)=
\begin{cases}
|k_*+y|^p-|k_*|^p,&k_n=k_*,\\
0,&|k_n|\to\infty.
\end{cases}
\]
In fact, for $|y|\le A$ and $|k_n|>2A$, the mean value theorem gives
$|V_n(y)|\le pA(|k_n|-A)^{p-1}\to0$ uniformly on this interval
as $n\to\infty$ in the case $|k_n|\to\infty$.

Since $(U_n)$ is bounded in the Hilbert space $\mathcal H_\rho$, there is a weakly convergent subsequence, still denoted by $(U_n)$, such that
$U_n\rightharpoonup U$ in $\mathcal H_\rho$ as $n\to\infty$. For every bounded spatial
interval $J$, $c_J:=\inf_J\rho>0$. Restriction to $(0,T)\times J$
is therefore a bounded map into $L^2((0,T)\times J)$, and the same subsequence converges
weakly there. For $\varphi(t,y)$ in $C_c^\infty((0,T)\times\R)$,
the bilinear distributional identity is
\[
\int_0^T\int_\R U_n(t,y)
 \bigl(-i\partial_t\varphi+\partial_y^2\varphi-V_n(y)\varphi\bigr)
 \,dy\,dt=0.
\]
Choose a bounded interval $J$ with $\supp\varphi\subset(0,T)\times J$.
The difference of the potential terms is
\[
\begin{aligned}
\int_0^T\int_\R (V_nU_n-V_\infty U)\varphi\,dy\,dt
 &=\int_0^T\int_J(V_n-V_\infty)U_n\varphi\,dy\,dt\\
 &\quad+\int_0^T\int_J V_\infty(U_n-U)\varphi\,dy\,dt .
\end{aligned}
\]
By $\norm{U_n}_{\mathcal H_\rho}=1$ and Cauchy-Schwarz inequality
\[
\left|\int_0^T\int_J(V_n-V_\infty)U_n\varphi\,dy\,dt\right|
 \leq c_J^{-1/2}\norm{V_n-V_\infty}_{L^\infty(J)}
       \norm\varphi_{L^2((0,T)\times J)}\longrightarrow0
       \quad\text{as }n\to\infty.
\]
The second integral tends to zero by weak convergence tested against
$\overline{V_\infty\varphi}\in L^2((0,T)\times J)$. 

The same weak convergence applies to the derivatives of $\varphi$. Thus
\[
i\partial_tU=-\partial_y^2 U+V_\infty U
\quad\text{in }\mathcal D'((0,T)\times\R).
\]

We apply Fact~\ref{fact-fc} to $U_n$ with the weight $\rho$ and the
cutoffs $\chi,\psi$ chosen above. Both $\rho$ and $\rho^{-1}$ are bounded
on bounded intervals, $\norm{U_n}_{\mathcal H_\rho}=1$, and the
time-independent potentials satisfy $|V_n(y)|\le |y|^p$.
After passing to a subsequence, we obtain the filtered limit:
\[
Z_n=\psi(D_t)(\chi U_n)\longrightarrow Z=\psi(D_t)(\chi U)
\quad\text{in }L^2_{\loc}(\R^2)\quad\text{as }n\to\infty.
\]
For every bounded interval $J$, the same fact gives
\[
e_J(S):=\sup_n\int_{|t|>S}\int_J|Z_n|^2\,dy\,dt\longrightarrow0
\qquad\text{as }S\to\infty.
\]
Passing to a further subsequence, we have $Z_n\to Z$ almost everywhere as $n\to\infty$.
For each $S>0$, Fatou's lemma gives
\[
\int_{|t|>S}\int_J|Z|^2\,dy\,dt
\leq\liminf_{n\to\infty}\int_{|t|>S}\int_J|Z_n|^2\,dy\,dt
\leq e_J(S).
\]
On $[-S,S]\times J$ we have strong $L^2$ convergence.
On $\{|t|>S\}\times J$ we use $|Z_n-Z|^2\leq2|Z_n|^2+2|Z|^2$. Hence
\[
\limsup_{n\to\infty}\norm{Z_n-Z}_{L^2(\R\times J)}^2\leq4e_J(S).
\]
Letting $S\to\infty$ after taking the upper limit in $n$ gives
$Z_n\to Z$ in $L^2(\R\times J)$ as $n\to\infty$.
Taking $J=Q_0$ in this conclusion and using \eqref{equ-nm} gives
\begin{equation}\label{equ-lm}
\norm Z_{L^2(\R\times Q_0)}^2
=\lim_{n\to\infty}\norm{Z_n}_{L^2(\R\times Q_0)}^2\geq\alpha>0.
\end{equation}

\Needspace{4\baselineskip}\medskip\noindent\textbf{Step 4. Vanishing of the filtered limit on a spatial set of positive measure.}\par\nopagebreak[4]
The functions $\ind_{E_n'}(y)$ belong to the closed unit ball of
$L^\infty(\R)$. By Banach--Alaoglu and the separability of $L^1(\R)$,
a further subsequence converges weak-star to a function $\zeta(y)$ in
$L^\infty(\R)$ as $n\to\infty$. For each nonnegative $h(y)$ in $L^1(\R)$,
\[
0\leq\int_\R\zeta(y)h(y)\,dy
=\lim_{n\to\infty}\int_\R\ind_{E_n'}(y)h(y)\,dy\leq\int_\R h(y)\,dy .
\]
Thus $0\leq\zeta\leq1$ almost everywhere. Testing against $\ind_{[x,x+L]}$ gives
\[
\int_x^{x+L}\zeta(y)\,dy=\lim_{n\to\infty}|E_n'\cap[x,x+L]|\geq\gamma L
\qquad(x\in\R).
\]
Let
$S=\{y\in(0,L):\zeta(y)>0\}$. Then
$\gamma L\leq\int_0^L\zeta(y)\,dy=\int_S\zeta(y)\,dy\leq|S|$.

For bounded intervals $K\subset\R_t$ and $J\subset\R_y$, by
\eqref{equ-no} and Plancherel,
\[
\begin{aligned}
\int_K\int_J\ind_{E_n'}(y)|Z_n(t,y)|^2\,dy\,dt
 &\leq\norm{\psi(D_t)(\ind_{E_n'\cap J}\chi U_n)}_{L^2(\R^2)}^2\\
 &\leq\norm\psi_\infty^2\norm\chi_\infty^2 c_J^{-1}\epsilon_n
 \longrightarrow0\qquad\text{as }n\to\infty.
\end{aligned}
\]
Let $h(y)=\ind_J(y)\int_K|Z(t,y)|^2\,dt$. Then $\norm h_{L^1(\R)}=\norm Z_{L^2(K\times J)}^2<\infty$. We have
\[
\begin{aligned}
&\left|\int_K\int_J\ind_{E_n'}|Z_n|^2\,dy\,dt-
         \int_K\int_J\zeta|Z|^2\,dy\,dt\right|\\
&=\left|\int_K\int_J\ind_{E_n'}(|Z_n|^2-|Z|^2)\,dy\,dt
 +\int_\R(\ind_{E_n'}-\zeta)h\,dy\right|\\
&\leq \norm{Z_n-Z}_{L^2(K\times J)}
   \bigl(\norm{Z_n}_{L^2(K\times J)}+\norm Z_{L^2(K\times J)}\bigr)
   +\left|\int_\R(\ind_{E_n'}-\zeta)h\,dy\right|
 \longrightarrow0\qquad\text{as }n\to\infty.
\end{aligned}
\]
The first bound on the last line follows from Cauchy--Schwarz. Its first term tends to
zero by strong convergence on $K\times J$, and its second by weak-star
convergence against the fixed $h$. Hence
$\int_K\int_J\zeta(y)|Z(t,y)|^2\,dy\,dt=0$.
Since the integrand is nonnegative and $K$ and $J$ are arbitrary, we obtain
$\zeta(y)|Z(t,y)|^2=0$ almost everywhere in $\R^2$.

\Needspace{4\baselineskip}\medskip\noindent\textbf{Step 5. Removal of the time filter and unique continuation.}\par\nopagebreak[4]
Since $S\subset(0,L)$ and $c_{(0,L)}>0$,
\[
\norm{\chi U}_{L^2(\R\times(0,L))}
 \leq \norm\chi_\infty c_{(0,L)}^{-1/2}\norm U_{\mathcal H_\rho},
\qquad
\norm{Z}_{L^2(\R\times(0,L))}
 \leq \norm\psi_\infty\norm{\chi U}_{L^2(\R\times(0,L))}.
\]
By Fubini, for almost every $y\in S$ the function $\chi U(\cdot,y)$
belongs to $L^2(\R_t)$ and is supported in $\supp\chi$. For almost
every such $y$, $Z(t,y)=0$ for almost every $t\in\R$, since
$\zeta(y)>0$ on $S$. For these $y$, up to a further null set, the
convolution representation gives
\[
\psi(D_t)(\chi U(\cdot,y))=K_\psi*(\chi U(\cdot,y))=Z(\cdot,y).
\]
The multiplier $\psi$ equals one on a nonempty open interval.
Fact~\ref{fact-fu}(i) gives $\chi(t)U(t,y)=0$ for
almost every $t$ and almost every $y\in S$. Since $\chi(t)>0$ on
$I_0$, $U=0$ almost everywhere on $I_0\times S$.

For every bounded interval $J$,
$\int_{I_0}\int_J|U(t,y)|^2\,dy\,dt\leq c_J^{-1}\norm U_{\mathcal H_\rho}^2$.
Thus $U\in L^2_{\loc}(I_0\times\R)$, and Step~3 gives the
distributional equation there with real, time-independent
$V_\infty\in L^\infty_{\loc}(\R)$. Since $I_0$ is connected and
$|S|>0$, Fact~\ref{fact-fu}(ii) yields $U=0$ almost everywhere on
$I_0\times\R$. Outside $I_0$ the cutoff is zero. Consequently
\[
\norm Z_{L^2(\R\times Q_0)}
\leq\norm\psi_\infty\norm{\chi U}_{L^2(\R\times Q_0)}=0,
\]
which contradicts the lower bound $\norm Z_{L^2(\R\times Q_0)}^2\geq
\alpha>0$ in \eqref{equ-lm}. Thus \eqref{equ-obs} holds with a constant depending only on $p,T,L,\gamma$.

\subsection{Necessity of thickness}\label{sec2-nec}
Suppose \eqref{equ-obs} holds with a finite constant for some
$T>0$, and assume for contradiction that $E$ is not thick. If $|E|=0$,
this is clearly a contradiction. Suppose $|E|>0$. By \cite[Lemma~7.1]{CLL}, there are centers $x_n$ with $|x_n|\to\infty$ and
\begin{equation}\label{equ-ti}
|E\cap(x_n-n,x_n+n)|<n^{-1}.
\end{equation}

Set $g(y)=\pi^{-1/4}e^{-y^2/2}$.
The self-adjoint translated and shifted operator
\[
H_n=-\partial_y^2+|x_n+y|^p-|x_n|^p
\]
has potential $\delta_n(y)=|x_n+y|^p-|x_n|^p$ converging locally
uniformly to zero as $n\to\infty$ by the mean value estimate above. Use the free solution
\[
z(t,y)=\pi^{-1/4}(1+2it)^{-1/2}
\exp\!\left(-\frac{y^2}{2(1+2it)}\right),
\]
with $(1+2it)^{-1/2}=1$ at $t=0$. Then $z(0)=g$ and
$i\partial_tz=-\partial_y^2z$. For $0\le t\le T$,
\[
|z(t,y)|^2=\frac{\exp(-y^2/(1+4t^2))}{\sqrt{\pi(1+4t^2)}}
\le C_Te^{-c_Ty^2},\qquad C_T=\pi^{-1/2},\quad c_T=(1+4T^2)^{-1}.
\]
For each integer $b\geq0$, $\partial_y^b z(t,y)$ is a
polynomial in $y$ times $z(t,y)$, with coefficients bounded
on $[0,T]$. Thus for all integers $a,b\geq0$,
\[
\sup_{0\leq t\leq T}\sup_{y\in\R}
(1+|y|)^a|\partial_y^b z(t,y)|\leq C_{a,b,T}.
\]

Since $|\delta_n(y)|\le|y|^p$ and $\delta_n(y)\to0$ as $n\to\infty$
for every fixed $y\in\R$, we obtain
\[
\varepsilon_n^2:=
\sup_{0\leq t\leq T}\norm{\delta_n z(t,\cdot)}_2^2
\leq C_T\int_\R|\delta_n(y)|^2 e^{-c_Ty^2}\,dy\longrightarrow0
\quad\text{as }n\to\infty
\]
by dominated convergence with majorant $|y|^{2p} e^{-c_Ty^2}$
inside the integral.

Fix $n$. Spatial Schwartz decay gives $z(t,\cdot)\in H^1(\R)$
and $|x_n+y|^{p/2}z$, $|x_n+y|^p z$,
$\partial_y^2z\in L^2(\R_y)$. 

Since $(i\partial_t-H_n)z=-\delta_nz$ and
$z(0,\cdot)=g$, Duhamel's formula gives
\[
e^{-itH_n}g-z(t,\cdot)
=-i\int_0^t e^{-i(t-r)H_n}\delta_n z(r,\cdot)\,dr,\qquad
\sup_{0\leq t\leq T}\norm{e^{-itH_n}g-z(t,\cdot)}_2
\leq T\varepsilon_n .
\]
Here the estimate follows from unitarity. Splitting the integral
at $|y|=n$ and using \eqref{equ-ti} gives
\[
\begin{aligned}
\int_0^T\norm{\ind_{E-x_n}e^{-itH_n}g}_2^2\,dt
&\leq 2\int_0^T\int_{E-x_n}|z(t,y)|^2\,dy\,dt
 +2T\sup_{0\leq t\leq T}\norm{e^{-itH_n}g-z(t,\cdot)}_2^2\\
&\leq 2TC_T n^{-1}
 +2TC_T\int_{|y|\geq n}e^{-c_Ty^2}\,dy
 +2T^3\varepsilon_n^2\longrightarrow0\qquad\text{as }n\to\infty.
\end{aligned}
\]

Set $f_n(x)=g(x-x_n)$ and $(T_nh)(x)=h(x-x_n)$.
For $h\in D(H_n)$, in the translated coordinates we have
\[
(T_n^{-1}H_mT_nh)(y)=-h''(y)+|x_n+y|^p h(y).
\]
Thus $T_n^{-1}H_mT_n=H_n+|x_n|^p I$, and the spectral theorem gives
\[
(e^{-itH_m}f_n)(x)
=e^{-it|x_n|^p}(e^{-itH_n}g)(x-x_n).
\]
And changing variables gives
\[
\norm{f_n}_2=1,\qquad
\int_0^T\int_E|(e^{-itH_m}f_n)(x)|^2\,dx\,dt
=\int_0^T\int_{E-x_n}|(e^{-itH_n}g)(y)|^2\,dy\,dt\longrightarrow0
\quad\text{as }n\to\infty.
\]
This contradicts \eqref{equ-obs}.

\begin{remark}
For $0<m<1/2$, Section~\ref{sec2-suf} gives observability of thick
sets in every positive time, uniformly for fixed $m,T,L,\gamma$.
Section~\ref{sec2-nec} shows that observation in even
one positive time implies thickness. At $m=1/2$, the same characterization
is given by \cite[Theorem~1.2]{CLL}.
\end{remark}

\section{Sufficient conditions for \texorpdfstring{$1/2<m\leq1$}{1/2 < m <= 1}}
\subsection{Interval density and a local oscillatory estimate}
\begin{lemma}\label{lem-ds}
For $1/2<m\leq1$ and $0<a<b$,
\begin{equation}\label{equ-dl}
0\leq b\,d_m(E,b)-a\,d_m(E,a)\leq b-a.
\end{equation}
The function $a\mapsto a\,d_m(E,a)$ is $1$-Lipschitz.
The positive-density scales are upward closed;
every $b>\kappa_m(E)$ is a positive-density scale if $\kappa_m(E)<\infty$,
and $d_m(E,\kappa_m(E))=0$ if $0<\kappa_m(E)<\infty$.
If $|E\mathbin{\triangle}F|<\infty$, then
$d_m(E,a)=d_m(F,a)$ for every $a>0$, and their critical scales agree.
\end{lemma}
\begin{proof}
The nested intervals differ in length $2(b-a)|x|^s$, so
\[
\frac{|E\cap I_{m,a}(x)|}{2|x|^s}
\leq\frac{|E\cap I_{m,b}(x)|}{2|x|^s}
\leq\frac{|E\cap I_{m,a}(x)|}{2|x|^s}+b-a.
\]
Taking lower limits gives \eqref{equ-dl} and continuity
of $d_m(E,\cdot)$ on $(0,\infty)$.
The infimum definition and upward closure give positivity above the critical
scale; continuity gives  $d_m(E,\kappa_m(E))=0$.
Finally,
\[
\bigl||E\cap I_{m,a}(x)|-|F\cap I_{m,a}(x)|\bigr|
\leq |E\mathbin{\triangle}F|.
\]
After division by $2a|x|^s$, the bound tends to zero as $|x|\to\infty$ because $s>0$.
\end{proof}

\begin{lemma}\label{lem-wt}
For a measurable $E\subset\R$, the following are equivalent:
\begin{enumerate}[label=\textup{(\roman*)},nosep]
\item $E$ is weakly thick, meaning
$\displaystyle\liminf_{R\to\infty}|E\cap[-R,R]|/(2R)>0$.
\item $d_1(E,a)>0$ for some $a>0$.
\item $d_1(E,b)>0$ for every $b>1$.
\end{enumerate}
Consequently $\kappa_1(E)\in[0,1]\cup\{+\infty\}$.
\end{lemma}
\begin{proof}
We proceed in the order 
(iii)$\implies$ (ii)$\implies$(i) $\implies$(iii).

Obviously, (iii) implies (ii). Then assume (ii). Choose $\gamma,x_0>0$ such that
$|E\cap I_{1,a}(x)|\geq2\gamma ax$ for $x\geq x_0$.
For $R\geq R_0:=(1+a)x_0$, set $x=R/(1+a)>0$.
Then
$I_{1,a}(x)=[(1-a)x,(1+a)x]\subset[-R,R]$, so
\begin{equation}\label{equ-wl}
|E\cap[-R,R]|\geq2\beta R\quad(R\geq R_0),\qquad
\beta=\frac{\gamma a}{1+a}>0.
\end{equation}
This holds for every sufficiently large $R$ and proves (i).

Now assume (i), with $\beta,R_0>0$ as in \eqref{equ-wl}.
For any $b>1$,
\[
[-(b-1)|x|,(b-1)|x|]\subset[x-b|x|,x+b|x|]=I_{1,b}(x).
\]
For $|x|\geq R_0/(b-1)$,  $R=(b-1)|x|\geq R_0$:
\[
\frac{|E\cap I_{1,b}(x)|}{2b|x|}
\geq\frac{|E\cap[-(b-1)|x|,(b-1)|x|]|}{2b|x|}
\geq\beta\frac{b-1}b.
\]
Hence $d_1(E,b)\geq\beta(b-1)/b>0$, proving (iii).
If a positive-density scale exists, their infimum is at most one;
otherwise it is infinite.
\end{proof}

\begin{lemma}\label{lem-ip}
Let $1/2<m\leq1$ and $d_m(E,a)>0$.
Allow every $a>0$ when $m<1$, and assume $0<a\leq1/4$ when $m=1$.
Then there is a countable family of bounded intervals
\[
\mathcal J=\{J_j:j\in\Z\}.
\]
More precisely, the intervals
are pairwise disjoint except possibly for common endpoints and
\[
\bigcup_{J\in\mathcal J}J=\R.
\]
There are constants $\gamma_0>0$ and $\ell_0>0$ such that every
$J\in\mathcal J$ satisfies
\begin{equation}\label{equ-pd}
|E\cap J|\geq\gamma_0|J|,\qquad |J|\geq\ell_0.
\end{equation}
For every noncentral interval $J=J_j$ with $j\neq0$,
\begin{equation}\label{equ-ps}
a|x|^s\leq |J|\leq4a|x|^s
\qquad(x\in J).
\end{equation}
The constants $\gamma_0,\ell_0$ and the central interval $J_0$
may depend on $m$, $a$, and $E$.
\end{lemma}
\begin{proof}
Choose $\gamma_a$ and $R_a$ as in \eqref{equ-dd}, so that
\begin{equation}\label{equ-pdd}
|E\cap I_{m,a}(x)|\geq2\gamma_a a|x|^s
\qquad(|x|\geq R_a).
\end{equation}
When $s<1$, enlarge $R_a$ so that
$ac^{s-1}\leq1/4$ for every $c\geq R_a$. When $s=1$, this
inequality holds because $a\leq1/4$. Thus the function
\[
F(c)=c-ac^s,\qquad c\geq R_a,
\]
is strictly increasing, since $F'(c)=1-as c^{s-1}\geq3/4$,
and $F(c)\to\infty$ as $c\to\infty$.

Choose $c_0$ so that $F(c_0)=R_a$. For $j\ge1$, given $b_j$,
choose $c_j$ so that
\begin{equation}\label{equ-pr}
F(c_j)=c_j-ac_j^s=b_j.
\end{equation}
In both cases set
\begin{equation}\label{equ-pre}
b_{j+1}=c_j+ac_j^s.
\end{equation}
Since $F$ is strictly increasing, tends to infinity and satisfies
$F(c)<c$, these equations have unique solutions with $c_0>R_a$ and
$c_j>b_j$ for $j\ge1$. In particular,
\[
b_{j+1}-b_j=2ac_j^s\geq2a R_a^s>0\qquad(j\ge1),
\]
so $b_j\to+\infty$ as $j\to\infty$.

Define the positive-side  intervals by
\[
J_j:=[b_j,b_{j+1}]
=[c_j-ac_j^s,c_j+ac_j^s]
=I_{m,a}(c_j),\qquad j\geq1.
\]
They satisfy
\[
[b_1,\infty)=\bigcup_{j\geq1}J_j,
\]
and consecutive intervals $J_j$ and $J_{j+1}$ meet only at $b_{j+1}$.
Since $c_j\geq R_a$, \eqref{equ-pdd} gives
\[
|E\cap J_j|\geq2\gamma_a ac_j^s
=\gamma_a|J_j|.
\]
For $x\in J_j$, the choice of $R_a$ gives
\[
\frac34c_j\leq x\leq\frac54c_j.
\]
As $0<s\leq1$, it follows that
\[
a|x|^s\leq a(5/4)^sc_j^s\leq2ac_j^s=|J_j|
\]
and
\[
|J_j|=2ac_j^s
\leq2a(4/3)^s x^s\leq4ax^s.
\]
This proves \eqref{equ-ps}.

For the negative side, reflect these intervals and define
\[
J_{-j}:=-J_j=[-b_{j+1},-b_j],\qquad j\geq1.
\]

Since $J_{-j}=I_{m,a}(-c_j)$, \eqref{equ-pdd} at $-c_j$
gives $|E\cap J_{-j}|\geq\gamma_a|J_{-j}|$. Reflection gives
\eqref{equ-ps} on $J_{-j}$.

We now define the single central interval
\[
J_0=[-b_1,b_1].
\]

By construction, $J_0$ joins the two half-lines and
\[
\R=\bigcup_{j\in\Z}J_j.
\]
Finally, $[R_a,b_1]=I_{m,a}(c_0)\subset J_0$. By \eqref{equ-pdd},
\[
|E\cap J_0|\geq|E\cap[R_a,b_1]|
\geq\gamma_a(b_1-R_a)\geq2\gamma_a a R_a^s>0.
\]
Thus the central interval also has positive observed proportion.
Set
\[
\gamma_0=\min\left\{\gamma_a,
\frac{|E\cap J_0|}{|J_0|}\right\},
\qquad
\ell_0=\min\{2a R_a^s,|J_0|\}.
\]
Every noncentral interval has length
$2ac_j^s\geq2a R_a^s$, while $J_0$ has length $|J_0|$.
This proves \eqref{equ-pd} and completes the proof.
\end{proof}

We next establish the local observation estimate used on these intervals.
\begin{lemma}\label{lem-lo}
Let $I$ be an interval of length $L>0$, let $u(x)\in H^2(I)$, and let
$S\subset I$ be measurable with $|S|\ge\gamma L$, where $0<\gamma\le1$.
Suppose $q(x)\in W^{1,\infty}(I;\R)$ satisfies
\begin{equation}\label{equ-oh}
k^2\le q(x)\le Ak^2,\qquad
\int_I\frac{|q'(x)|}{q(x)}\,dx\le M_0,\qquad
\frac{|q'(x)|}{q(x)^{3/2}}\le1\quad\text{a.e. on }I
\end{equation}
for some $k>0$, $A\ge1$ and $M_0\ge0$. There exist constants $K_0=K_0(A)>0$
and $C_{\rm obs}=C_{\rm obs}(A,M_0)>0$ such that, whenever $\gamma kL\ge K_0$,
\begin{equation}\label{equ-lo}
\norm u_{L^2(I)}^2\le C_{\rm obs}\gamma^{-3}\norm u_{L^2(S)}^2+
4e^{M_0}\frac{L^2}{k^2}\norm{u''+qu}_{L^2(I)}^2.
\end{equation}
\end{lemma}
\begin{proof}
First suppose $u$ is real and set $f(x)=u''(x)+q(x)u(x)$.
Use the $C^1$ representative of $u$ on $\overline I$, whose derivative
is absolutely continuous, and the Lipschitz representative of $q$. Put
\[
v(x)=\frac{u'(x)}{\sqrt{q(x)}},\qquad
R_u(x)=\sqrt{u(x)^2+v(x)^2}.
\]
Since $q\ge k^2$, the vector $(u,v)$ and its length $R_u$ are absolutely
continuous on $\overline I$. Choose $x_*\in\overline I$ such that
$r_*:=R_u(x_*)=\min_{\overline I}R_u$. Then
\[
(R_u^2)'=2uu'+\frac{2u'u''}{q}-\frac{q'(u')^2}{q^2}
=\frac{2u'f}{q}-\frac{q'(u')^2}{q^2}
\qquad\text{a.e. on }I.
\]
At each interior zero where $R_u$ is differentiable, $R_u'=0$ since $R_u\ge0$.
Where $R_u>0$, dividing by $2R_u$ and using $|u'|\le\sqrt q\,R_u$ gives
\[
|R_u'|\le\frac{|q'|}{2q}R_u+\frac{|f|}{k}
\qquad\text{a.e. on }I.
\]
For $x\in\overline I$, let $J_x$ be the interval between $x_*$ and $x$.
Since $k^{-1}\int_{J_x}|f(s)|\,ds\le D_f:=k^{-1}\sqrt L\,\norm f_{L^2(I)}$,
integration along $J_x$ and Gronwall's inequality give
\[
R_u(x)\le(r_*+D_f)\exp\!\left(\int_{J_x}\frac{|q'(s)|}{2q(s)}\,ds\right)
\le e^{M_0/2}(r_*+D_f).
\]
\smallskip\noindent\textit{Case I. $r_*\le D_f$.} Then
\[
\norm u_{L^2(I)}^2\le L\max_{\overline I}R_u^2
\le 4e^{M_0}\frac{L^2}{k^2}\norm f_{L^2(I)}^2.
\]

\smallskip\noindent\textit{Case II. $r_*>D_f$.} Since $R_u(x)\ge r_*>0$, the normalized
vector $(v,u)/R_u$ traces an absolutely continuous curve on the unit circle.
Hence there exists an absolutely continuous angle function \(\theta\colon I\to\mathbb{R}\) such that
\[
v(x)=R_u(x)\cos\theta(x),\qquad u(x)=R_u(x)\sin\theta(x).
\]
Differentiating gives $vu'-uv'=R_u^2\theta'$ almost everywhere. Since
$u'=\sqrt q\,v$ and $v'=-\sqrt q\,u-q'v/(2q)+f/\sqrt q$, it follows that
\begin{equation}\label{equ-pa}
\theta'=\frac{vu'-uv'}{R_u^2}
=\sqrt q+\frac{q'}{4q}\sin(2\theta)-\frac{f}{\sqrt q\,R_u}\sin\theta.
\end{equation}
Set $\theta'=b+e$, where $b(x)=\sqrt q+(q'/4q)\sin(2\theta)$
and $e(x)=-f\sin\theta/(\sqrt q\,R_u)$.
The bounds on $q$ and $q'$ give
\[
\frac{3k}4\le b(x)\le\frac{5\sqrt A\,k}4,\qquad
\int_I|e(x)|\,dx\le\frac1{kr_*}\int_I|f(x)|\,dx
\le\frac{D_f}{r_*}<1.
\]

For $0<\epsilon<1/4$, set $B_\epsilon=\{x\in I:|\sin\theta(x)|<\epsilon\}$.
Fix a smooth function $\zeta(r)$ with values in $[0,1]$, equal to one on
$(-\infty,1]$ and zero on $[4,\infty)$. Define
\[
G_\epsilon(z)=\zeta(\sin^2 z/\epsilon^2),\quad
\mu_\epsilon=\frac1\pi\int_0^\pi G_\epsilon(z)\,dz,\quad
P_\epsilon(z)=\int_0^z(G_\epsilon(r)-\mu_\epsilon)\,dr.
\]
The function $G_\epsilon$ is $\pi$-periodic, equals one for
$|\sin z|\le\epsilon$ and is zero for $|\sin z|\ge2\epsilon$.
Hence $0\le\pi\mu_\epsilon\le2\arcsin(2\epsilon)\le2\pi\epsilon$.
The integrand defining $P_\epsilon$ has mean zero over a period, so
$P_\epsilon$ is $\pi$-periodic and
$\norm{P_\epsilon}_{L^\infty(\R)}\le2\pi\mu_\epsilon$.
Let $\overline I=[\alpha,\beta]$. Then
$|\theta(\beta)-\theta(\alpha)|\le 5\sqrt A\,kL/4+1$.
The chain rule therefore gives
\[
\begin{aligned}
\left|\int_I G_\epsilon(\theta(x))\theta'(x)\,dx\right|
&=\left|\mu_\epsilon(\theta(\beta)-\theta(\alpha))
+P_\epsilon(\theta(\beta))-P_\epsilon(\theta(\alpha))\right|\\
&\le 2\epsilon(5\sqrt A\,kL/4+1)+8\pi\epsilon
\le C_A\epsilon(kL+1),
\end{aligned}
\]
where $C_A=32\sqrt A$. Using $b=\theta'-e$ and $0\le G_\epsilon\le1$, we obtain
\[
\begin{aligned}
\frac{3k}{4}|B_\epsilon|
&\le \frac{3k}{4}\int_I G_\epsilon(\theta(x))\,dx
\le \int_I G_\epsilon(\theta(x))b(x)\,dx\\
&\le C_A\epsilon(kL+1)+\int_I|e(x)|\,dx
\le C_A\epsilon(kL+1)+1.
\end{aligned}
\]
Dividing by $3k/4$ gives $|B_\epsilon|\le C_{\rm ph}(\epsilon L+1/k)$
with $C_{\rm ph}=48\sqrt A$. Set $c_1=(4C_{\rm ph})^{-1}$,
$\epsilon=c_1\gamma$ and $K_0=4C_{\rm ph}$. Since $\gamma kL\ge K_0$,
we have $C_{\rm ph}/k\le\gamma L/4$. Thus
\[
|B_\epsilon|\le C_{\rm ph}c_1\gamma L+C_{\rm ph}/k\le\gamma L/2,\qquad
|S\setminus B_\epsilon|\ge |S|-|B_\epsilon|\ge\gamma L/2.
\]
For $x\in S\setminus B_\epsilon$ we have
$|u(x)|=R_u(x)|\sin\theta(x)|\ge r_*c_1\gamma$. Consequently,
\[
\int_S|u(x)|^2\,dx\ge (r_*c_1\gamma)^2|S\setminus B_\epsilon|
\ge (c_1^2/2)\gamma^3 Lr_*^2.
\]
Since $D_f<r_*$, we have $R_u(x)\le 2e^{M_0/2}r_*$ on $\overline I$. Hence
\[
\norm u_{L^2(I)}^2\le 4e^{M_0}Lr_*^2
\le \frac{8e^{M_0}}{c_1^2}\gamma^{-3}\norm u_{L^2(S)}^2.
\]
Combining the two cases gives \eqref{equ-lo} with
$C_{\rm obs}=8e^{M_0}/c_1^2$, where $c_1$ depends only on $A$.
Since $q$ is real, applying the estimate to the real and imaginary parts of $u$
and adding proves \eqref{equ-lo}.
\end{proof}

\subsection{Spectral estimates and a single-scale resolvent estimate}
Throughout this section let $1/2<m\leq1$.
The symbol $C_m$ denotes a positive constant depending only on $m$,
whose value may change from line to line.
The final forcing constant is denoted by $A_m$.

For the harmonic oscillator, Koch--Tataru \cite{KT} give $L^p$ eigenfunction
bounds. Spectral inequalities for finite Hermite expansions are studied in
\cite{BJPS,MPSherm,DSVhermite}, and for Schrödinger operators with power growth
potentials in \cite{DSVpower}. Egidi--Seelmann~\cite{ES} give an abstract
Logvinenko--Sereda theorem for spectral subspaces under Bernstein-type
assumptions. Here we need a bound on the mass near and beyond
the turning points for eigenfunctions in a short spectral interval.

The eigenvalues of $H_m=-\partial_x^2+|x|^{2m}$ on $L^2(\R)$ are simple
and satisfy $0<\lambda_1<\lambda_2<\cdots$ \cite[Section~7.1]{DM1}.
Let $\varphi_n$ be a real eigenfunction with
$H_m\varphi_n=\lambda_n\varphi_n$ and $\norm{\varphi_n}_2=1$.
Set $r_n=\lambda_n^{1/(2m)}$. The even convex
potential $V(x)=|x|^{2m}$ satisfies $xV'(x)=2mV(x)$ for $x\ne0$,
so $V\in\mathcal P_1(2m)$ in \cite[Definitions~6.1 and 7.5]{DM1}
and $V\in\mathcal P_{1,\mathrm{cv}}(2m)$ in \cite[Definition~3.3]{DM2}.

\medskip\noindent\textbf{Spectral fact 1. Eigenvalues and spectral gaps.}
By \cite[Theorem~7.6]{DM1},
\begin{equation}\label{equ-si}
\sqrt{\lambda_n}\,|\{V\leq\lambda_n\}|\asymp_m n,\qquad
\lambda_j-\lambda_i\asymp_m\frac{\lambda_j}j(j-i)
\quad(1\leq i<j).
\end{equation}
Here $X\asymp_mY$ means two-sided inequalities with positive constants
depending only on $m$.
Since $|\{V\leq\lambda_n\}|=2r_n$, the first estimate gives
$n\asymp_mr_n^{m+1}$.
The second with $j=n+1$ gives
$(\lambda_{n+1}-\lambda_n)/\lambda_{n+1}=O_m(n^{-1})$,
so $\lambda_{n+1}/\lambda_n\to1$ as $n\to\infty$. Therefore
\begin{equation}\label{equ-sg}
n\asymp_mr_n^{m+1},\qquad
\lambda_{n+1}-\lambda_n\asymp_mr_n^{m-1}.
\end{equation}

\medskip\noindent\textbf{Spectral fact 2. Pointwise eigenfunction estimates.}
For $V\in\mathcal P_{1,\mathrm{cv}}(2m)$,
\cite[Theorem~3.5(iii)]{DM2} gives, with $\alpha=1/4$, for every $n\geq1$ and $x\in\R$,
\[
|\varphi_n(x)|\leq\frac{C_m}{|\{V\leq\lambda_n\}|^{1/2}}
\min\{n^{1/6},|1-V(x)/\lambda_n|^{-1/4}\}.
\]
Here $n^{2\alpha/3}=n^{1/6}$. Squaring the estimate and using
$|\{V\leq\lambda_n\}|=2r_n$ and
$n\asymp_m r_n^{m+1}$ by \eqref{equ-sg}, we obtain, for $|x|\leq r_n$,
\begin{equation}\label{equ-ef}
|\varphi_n(x)|^2\leq\frac{C_m}{r_n}
\min\left\{r_n^{(m+1)/3},
\left(1-\frac{|x|^{2m}}{\lambda_n}\right)^{-1/2}\right\}.
\end{equation}
At $x=\pm r_n$, the second entry in the minimum is understood as
$+\infty$, and the first gives
\[
|\varphi_n(\pm r_n)|^2\leq C_m r_n^{(m-2)/3}.
\]
For $|x|<r_n$, the second entry instead gives
\[
|\varphi_n(x)|^2\leq C_m r_n^{m-1}
\bigl(\lambda_n-V(x)\bigr)^{-1/2}.
\]

\begin{lemma}\label{lem-fm}
For all sufficiently large $n$,
\begin{equation}\label{equ-fm}
\int_{|x|>r_n}|\varphi_n(x)|^2\,dx\leq C_mr_n^{-(m+1)/3}.
\end{equation}
\end{lemma}
\begin{proof}
Set $r=r_n$, $\lambda=\lambda_n$ and $\varphi=\varphi_n$.
On $(r,\infty)$ the potential $V-\lambda$ is positive.
Since $\varphi\in L^2(r,\infty)$, \cite[Proposition~5.2]{DM1} applies to
$-\varphi''+(V-\lambda)\varphi=0$. Choose the sign of $\varphi$ so that
$\varphi>0$ and $\varphi'<0$ on $(r,\infty)$. In particular, $\varphi(r)>0$.

Set $G=(\varphi')^2-(V-\lambda)\varphi^2$. The eigenfunction equation gives
$G'=-V'\varphi^2\leq0$ on $(r,\infty)$. Since $\varphi\in Q_m$,
\[
|G|\leq |\varphi'|^2+V|\varphi|^2\in L^1(r,\infty).
\]
Thus $G$ is nonincreasing and integrable, so $G(x)\to0$ as $x\to\infty$.
Integrating $G'$ gives
\[
G(x)=\int_x^\infty V'(t)\varphi(t)^2\,dt\geq0,\qquad
-\frac{\varphi'(x)}{\varphi(x)}\geq\sqrt{V(x)-\lambda}
\quad(x>r).
\]
Since $V(r)=\lambda$ and $V$ is convex, for $y\ge0$ we have
$V(r+y)-\lambda\ge V'(r)y=2mr^s y$. Thus for $t\ge0$,
\[
\log\frac{\varphi(r)}{\varphi(r+t)}
=\int_0^t-\frac{\varphi'(r+y)}{\varphi(r+y)}\,dy
\ge\sqrt{2m}\,r^{s/2}\int_0^t y^{1/2}\,dy
=\frac23\sqrt{2m}\,r^{s/2}t^{3/2}.
\]
Exponentiating and squaring, we obtain by \eqref{equ-ef} at $r$
\[
\begin{aligned}
|\varphi(r+t)|^2
&\le |\varphi(r)|^2\exp\!\left(-\frac43\sqrt{2m}\,r^{s/2}t^{3/2}\right)\\
&\le C_mr^{(m-2)/3}\exp\!\left(-\frac43\sqrt{2m}\,r^{s/2}t^{3/2}\right).
\end{aligned}
\]
Consequently,
\[
\begin{aligned}
\int_r^\infty|\varphi(x)|^2\,dx
&=\int_0^\infty|\varphi(r+t)|^2\,dt\\
&\leq C_m r^{(m-2)/3}
 \int_0^\infty e^{-c_m r^{s/2}t^{3/2}}\,dt\\
&= C_m r^{(m-2-s)/3}\int_0^\infty e^{-c_m v^{3/2}}\,dv
\leq C_m r^{-(m+1)/3},
\end{aligned}
\]
where $t=r^{-s/3}v$. Since $V$ is even,
$\widetilde\varphi(x)=\varphi(-x)$ satisfies the same eigenfunction equation.
Applying the above estimate to $\widetilde\varphi$ gives
\[
\int_{-\infty}^{-r}|\varphi(x)|^2\,dx
=\int_r^\infty|\widetilde\varphi(x)|^2\,dx\leq C_m r^{-(m+1)/3}.
\]
Adding the estimates on the two half-lines proves \eqref{equ-fm}.
\end{proof}

For $\lambda>0$, $0<K\leq\lambda/4$ and $0<w\leq K/2$, put
\begin{equation}\label{equ-scd}
r=\lambda^{1/(2m)},\quad b=r^s,\quad
\Omega_{\lambda,K}=\{V\leq\lambda-K\},\quad
P_{\lambda,w}=\ind_{[\lambda-w,\lambda+w]}(H_m).
\end{equation}
Explicitly,
$P_{\lambda,w}u=\sum_{|\lambda_n-\lambda|\leq w}
\langle u,\varphi_n\rangle\varphi_n$.

\begin{lemma}\label{lem-sc}
For all sufficiently large $\lambda$, $0<K\leq\lambda/4$ and $0<w\leq K/2$,
\begin{equation}\label{equ-sc}
\norm{\ind_{\Omega_{\lambda,K}^c}P_{\lambda,w}}_{2\to2}^2
\leq C_m(1+wr^{1-m})
\bigl(r^{-m}\sqrt K+r^{-(m+1)/3}\bigr).
\end{equation}
\end{lemma}
\begin{proof}
Let $\mathcal I=\{n:|\lambda_n-\lambda|\le w\}$. If $\mathcal I$ is empty,
then $P_{\lambda,w}=0$ and the estimate is immediate.
For arbitrary $f\in L^2(\mathbb R)$, set
\[
P_{\lambda,w}f
=\sum_{n\in\mathcal I}c_n\varphi_n,
\qquad
c_n=\langle f,\varphi_n\rangle.
\]
By Bessel's inequality,
\[
\sum_{n\in\mathcal I}|c_n|^2
=\|P_{\lambda,w}f\|_2^2
\leq\|f\|_2^2.
\]
By the Cauchy--Schwarz inequality,
\[
|P_{\lambda,w}f(x)|^2
\leq
\left(\sum_{n\in\mathcal I}|c_n|^2\right)
\left(\sum_{n\in\mathcal I}|\varphi_n(x)|^2\right)
\leq
\|f\|_2^2\sum_{n\in\mathcal I}|\varphi_n(x)|^2.
\]
Integrating over $\Omega_{\lambda,K}^c$, we obtain
\begin{equation}\label{equ-cps}
\|\mathbf1_{\Omega_{\lambda,K}^c}P_{\lambda,w}f\|_2^2
\leq
\|f\|_2^2
\sum_{n\in\mathcal I}
\int_{\Omega_{\lambda,K}^c}|\varphi_n(x)|^2\,dx.
\end{equation}

Fix $n\in\mathcal I$. Since $\lambda_n\ge\lambda-w>\lambda-K$,
the integration region is the disjoint union of
\[
A_n=\{\lambda-K<V\leq\lambda_n\},
\qquad
F_n=\{V>\lambda_n\}=\{|x|>r_n\}.
\]
On $A_n$ we have $0\leq\lambda_n-V\leq K+w\leq3K/2$.
The bounds $3\lambda/4<V\leq\lambda+w\leq9\lambda/8$ imply
$|x|\asymp_m r$ and $|V'|\asymp_m r^s$.
Also $r_n\asymp_m r$ for all $n\in\mathcal I$.
By \eqref{equ-ef}, for every $x\in A_n$ with $V(x)<\lambda_n$,
\[
|\varphi_n(x)|^2
\leq C_mr^{m-1}(\lambda_n-V(x))^{-1/2}.
\]
Changing variables to $z=\lambda_n-V(x)$ on each half-line gives
\[
\begin{aligned}
\int_{A_n}|\varphi_n(x)|^2\,dx
&\leq C_mr^{m-1-s}
\int_0^{3K/2}z^{-1/2}\,dz\\
&\leq C_mr^{-m}\sqrt K.
\end{aligned}
\]
Lemma~\ref{lem-fm} gives
\[
\int_{F_n}|\varphi_n(x)|^2\,dx
\leq C_mr^{-(m+1)/3}.
\]
Consequently,
\begin{equation}\label{equ-com}
\int_{\Omega_{\lambda,K}^c}|\varphi_n(x)|^2\,dx
\leq C_m\left(r^{-m}\sqrt K+r^{-(m+1)/3}\right).
\end{equation}

If the spectral interval $[\lambda-w,\lambda+w]$ contains $N\geq2$
eigenvalues, let their consecutive indices be
$n_1<\cdots<n_N$. By \eqref{equ-sg},
\[
\begin{aligned}
2w
&\geq\lambda_{n_N}-\lambda_{n_1}\\
&=\sum_{j=1}^{N-1}
(\lambda_{n_{j+1}}-\lambda_{n_j})\\
&\geq c_m(N-1)r^{m-1}.
\end{aligned}
\]
Thus
\begin{equation}\label{equ-cwc}
N\leq1+C_mwr^{1-m},
\end{equation}
which also holds for $N=0,1$. By \eqref{equ-cps},
\eqref{equ-com} and \eqref{equ-cwc}, for every $f\in L^2(\R)$,
\[
\begin{aligned}
\|\mathbf1_{\Omega_{\lambda,K}^c}P_{\lambda,w}f\|_2^2
&\leq C_mN
\left(r^{-m}\sqrt K+r^{-(m+1)/3}\right)\|f\|_2^2\\
&\leq C_m(1+wr^{1-m})
\left(r^{-m}\sqrt K+r^{-(m+1)/3}\right)\|f\|_2^2.
\end{aligned}
\]
This proves \eqref{equ-sc}.
\end{proof}

We now combine the interval and spectral estimates to prove the following resolvent bound.
\begin{proposition}\label{prop-re}
Let $1/2<m\le1$. There exist $A_m>0$, depending only on $m$, and a
universal $a_0\in(0,1/4]$ such that the following holds.
Assume $a>0$ if $m<1$, and $0<a\le a_0$ if $m=1$.
If $d_m(E,a)>0$, there are finite constants $N_a$ and $\Lambda_a$,
which may depend on $m$, $a$, and $E$, such that
\begin{equation}\label{equ-re}
\norm u_2^2\leq N_a\norm{\ind_Eu}_2^2+
A_ma\norm{(H_m-\lambda)u}_2^2,\qquad
\lambda\geq\Lambda_a,\quad u\in D(H_m).
\end{equation}
\end{proposition}
\begin{proof}
Fix $a$ with $d_m(E,a)>0$ as above and take the partition of
Lemma~\ref{lem-ip}. For $u\in D(H_m)$ set $f=(H_m-\lambda)u$ and
\begin{equation}\label{equ-kc}
K=Dab^2,
\end{equation}
where $D=D(m)$ is fixed sufficiently large.

\Needspace{4\baselineskip}\medskip\noindent\textbf{Step 1. The resolvent estimate on the allowed region.}\par\nopagebreak[4]
Set $\Omega=\Omega_{\lambda,K}=\{V\le\lambda-K\}$. From the partition $\mathcal J$
of Lemma~\ref{lem-ip} select the family
\[
\mathcal J_\Omega=\{J\in\mathcal J:J\cap\Omega\ne\varnothing\}.
\]
For each selected $J$, choose $x_J\in J\cap\Omega$, so $|x_J|\leq r$.
For a noncentral interval, \eqref{equ-ps} gives $|J|\leq4ab$.
Set $G=4$. Take $r\ge1$ sufficiently large that
\[
 |J_0|\le ar^s,\qquad Gar^{s-1}\le1.
\]
The first inequality holds for large $r$ since $s>0$. The second holds for
large $r$ if $m<1$, and for $a\le1/G$ if $m=1$.
Thus $|J|\le Gab$ for every selected interval. For $x\in J$,
\[
 |x|\le |x_J|+|J|\le r+Gar^s\le2r.
\]
Fix $D\ge\max\{1,2(2m\,2^sG)\}$. Since $q=\lambda-V$ and
$|V'|\le 2m(2r)^s$ on $[-2r,2r]$,
\[
 \osc_J q=\osc_J V\le \sup_J|V'|\,|J|
 \le 2m\,2^sG\,ab^2\le K/2.
\]
Moreover, $q(x_J)\ge K$.
With $k_J^2=\inf_Jq$ this gives
\begin{equation}\label{equ-qc}
k_J^2\geq K/2,\qquad \sup_Jq\leq2\inf_Jq.
\end{equation}

Moreover, $q$ is monotone on each side of the origin. Hence
\[
\int_J\frac{|q'(x)|}{q(x)}\,dx
\leq 2\log\frac{\sup_Jq}{\inf_Jq}
\leq2\log2,
\]
where the factor $2$ accounts for the two parts of $J$.
The remaining hypotheses of Lemma~\ref{lem-lo} follow from
\[
\sup_J\frac{|q'|}{q^{3/2}}
\leq2m2^s(2/D)^{3/2}a^{-3/2}r^{-2s},\qquad
k_J|J|\geq\sqrt{Da/2}\,r^s\ell_0.
\]
Increase $r$ so that
\[
r^{2s}\geq2m2^s(2/D)^{3/2}a^{-3/2},\qquad
r^s\geq\sqrt{2/D}\,\frac{K_0}{\gamma_0\ell_0\sqrt a}.
\]
This is possible since $s>0$.
Since $u\in H^2_{\loc}$ and $f=-(u''+qu)$,
Lemma~\ref{lem-lo} applies on each selected interval with
$A=2$, $M_0=2\log2$ and observed set $E\cap J$.
Moreover,
\[
\frac{|J|^2}{k_J^2}\leq
\frac{G^2a^2b^2}{Dab^2/2}=\frac{2G^2}D\,a.
\]
Since $\Omega\subseteq\bigcup_{J\in\mathcal J_\Omega}J$,
$\int_\Omega|u|^2\le\sum_{J\in\mathcal J_\Omega}\int_J|u|^2$.
Summing \eqref{equ-lo} with
$M_0=2\log2$ and observation sets $E\cap J$ therefore gives
\begin{equation}\label{equ-ar}
\int_\Omega|u|^2
\leq N_a'\norm{\ind_Eu}_2^2+C_m^{\rm al}a\norm f_2^2,
\end{equation}
with
\[
N_a'=C_{\rm obs}(2,2\log2)\gamma_0^{-3},\qquad
C_m^{\rm al}=\frac{2(4e^{M_0})G^2}D=\frac{32G^2}D,
\qquad 4e^{M_0}=4e^{2\log2}=16.
\]
Here $\gamma_0>0$ is the lower bound for $|E\cap J|/|J|$ in Lemma~\ref{lem-ip}.

\Needspace{4\baselineskip}\medskip\noindent\textbf{Step 2. Control of the mass outside the allowed region.}\par\nopagebreak[4]
Choose $w=\epsilon a^{-1/2}$, with $\epsilon=\epsilon(m)>0$.
Then $w^{-2}=\epsilon^{-2}a$. Since $K=Dar^{2s}$ and $s>0$, we have
$w\le K/2$ for all sufficiently large $r$. Also $K/\lambda=Dar^{2m-2}$,
which tends to zero as $r=\lambda^{1/(2m)}\to\infty$ for fixed $D,a$
if $m<1$, and equals $Da$ if $m=1$.
Thus $K\le\lambda/4$ for large $r$, provided $a\le1/(4D)$ at $m=1$.
In \eqref{equ-sc},
\[
r^{-m}\sqrt K=\sqrt{Da}\,r^{m-1},\qquad
wr^{1-m}r^{-m}\sqrt K=\sqrt D\,\epsilon.
\]
Thus
\begin{equation}\label{equ-ce}
\begin{split}
\norm{\ind_{\Omega^c}P_{\lambda,w}}_{2\to2}^2
\leq C_m^{\rm col}\bigl[
&\sqrt D\,\epsilon+\sqrt{Da}\,r^{m-1}\\
&+r^{-(m+1)/3}+\epsilon a^{-1/2}r^{(2-4m)/3}\bigr].
\end{split}
\end{equation}
With $m,D,a,\epsilon$ fixed, $r=\lambda^{1/(2m)}\to\infty$
is equivalent to $\lambda\to\infty$. We have
\[
r^{-(m+1)/3}\to0,\qquad
\epsilon a^{-1/2}r^{(2-4m)/3}\to0
\quad\text{as }r\to\infty,
\]
since $2-4m<0$. If $1/2<m<1$, then
$\sqrt{Da}\,r^{m-1}\to0$ as $r\to\infty$ as well. At $m=1$, this term is
$\sqrt{Da}$ and is controlled by the choice of $a_0$.
Choose $0<\epsilon\leq(64C_m^{\rm col}\sqrt D)^{-1}$.
At $m=1$ choose in addition
\begin{equation}\label{equ-a0}
0<a_0\leq
\min\left\{\frac14,\frac1G,\frac1{4D},
\frac1{(64C_1^{\rm col})^2D}\right\}.
\end{equation}
With these choices, each nondecaying term on the right of
\eqref{equ-ce}, including the factor $C_m^{\rm col}$, is at most $1/64$.
Since the remaining terms tend to zero as $r\to\infty$, choose $r_0\ge1$ so that for
every $r\ge r_0$ they satisfy the same bound, all preceding large-$r$
conditions hold, and $\lambda=r^{2m}$ is in the range of Lemma~\ref{lem-sc}.
Taking $\Lambda_a=r_0^{2m}$, we obtain for every $\lambda\ge\Lambda_a$
\begin{equation}\label{equ-csml}
\norm{\ind_{\Omega^c}P_{\lambda,w}}_{2\to2}^2\le\tfrac1{16}.
\end{equation}

With $P=P_{\lambda,w}$ and $u=\sum_nc_n\varphi_n$, orthogonality gives
\[
\norm{(1-P)u}_2^2
=\sum_{|\lambda_n-\lambda|>w}|c_n|^2
\leq w^{-2}\sum_n|\lambda_n-\lambda|^2|c_n|^2
=\epsilon^{-2}a\norm f_2^2.
\]
Using $u=Pu+(1-P)u$, contraction of spatial restriction and \eqref{equ-csml},
\begin{equation}\label{equ-cm}
\norm{\ind_{\Omega^c}u}_2^2
\leq2\norm{\ind_{\Omega^c}P}_{2\to2}^2\norm u_2^2
+2\norm{(1-P)u}_2^2
\leq\frac18\norm u_2^2+2\epsilon^{-2}a\norm f_2^2.
\end{equation}
Add \eqref{equ-ar} and \eqref{equ-cm}, then absorb
$\norm u_2^2/8$. This proves \eqref{equ-re} with
\[
N_a=\frac87N_a',\qquad
A_m=\frac87\left(\frac{32G^2}D+2\epsilon^{-2}\right).
\]
\end{proof}

\subsection{From high energy to observation in time}
We adapt the time-Fourier argument of Miller \cite[Theorem~5.1]{Miller}
to a resolvent estimate assumed only at high energies.
For spectral criteria for exact observability and the operator-semigroup
framework, see \cite{RTTT,TW}.
The low-energy remainder is removed by finite-dimensional compactness
and injectivity of the observation operator on each eigenspace.

\begin{proposition}\label{prop-suf}
Let $1/2<m\leq1$ and let $a$ be allowed in Proposition~\ref{prop-re}.
If $d_m(E,a)>0$, then \eqref{equ-obs} holds for every
\begin{equation}\label{equ-sft}
T>\pi\sqrt{A_ma}.
\end{equation}
\end{proposition}
\begin{proof}
Fix $a$ as in Proposition~\ref{prop-re} with $d_m(E,a)>0$ and take $T>\pi\sqrt{A_ma}$.

\Needspace{4\baselineskip}\medskip\noindent\textbf{Step 1. A time-observation estimate with a finite-rank remainder.}\par\nopagebreak[4]
Choose $R_0>\max\{0,\Lambda_a\}+(A_ma)^{-1/2}$ and $P_0=\ind_{[0,R_0]}(H_m)$.
For $\lambda<\Lambda_a$, the spectral theorem gives
$\norm{(1-P_0)v}_2^2\leq A_ma\norm{(H_m-\lambda)v}_2^2$.
By \eqref{equ-re} for $\lambda\ge\Lambda_a$ and orthogonality for $\lambda<\Lambda_a$,
we have, for every $\lambda\in\R$ and $v\in D(H_m)$,
\begin{equation}\label{equ-abr}
\norm v_2^2\leq N_a\norm{\ind_Ev}_2^2+
A_ma\norm{(H_m-\lambda)v}_2^2+\norm{P_0v}_2^2.
\end{equation}

Let $\chi_T(t)=\sin(\pi t/T)$ on $(0,T)$ and zero elsewhere.
For $f$ in a bounded spectral subspace, put $v(t)=\chi_T(t)U_m(t)f$.
For almost everywhere $\tau$, our Fourier convention gives
\[
(H_m+\tau)\widehat v(\tau)=-i\,\widehat{\chi_T'U_m(t)f}(\tau).
\]
By Plancherel's theorem,
\[
\int_\R\norm{(H_m+\tau)\widehat v(\tau)}_2^2\,d\tau
=\int_0^T|\chi_T'(t)|^2\norm{U_m(t)f}_2^2\,dt
=\frac{\pi^2}{2T}\norm f_2^2.
\]
Apply \eqref{equ-abr} to $\widehat v(\tau)$ with $\lambda=-\tau$ and integrate in $\tau$.
Since $\int_0^T\chi_T^2=T/2$, $|\chi_T|\leq1$ and $P_0U_m(t)=U_m(t)P_0$,
\begin{equation}\label{equ-tr}
\left(\frac T2-\frac{A_ma\pi^2}{2T}\right)\norm f_2^2
\leq N_a\int_0^T\norm{\ind_E U_m(t)f}_2^2\,dt+\frac T2\norm{P_0f}_2^2.
\end{equation}
Since $\int_0^T\norm{\ind_E U_m(t)h}_2^2\,dt\le T\norm h_2^2$ and $P_0$ is bounded,
spectral truncation extends \eqref{equ-tr} to every $f\in L^2(\R)$.
By \eqref{equ-sft}, $\alpha_T:=T/2-A_ma\pi^2/(2T)>0$.

\Needspace{4\baselineskip}\medskip\noindent\textbf{Step 2. Removal of the remainder by compactness and uniqueness.}\par\nopagebreak[4]
If \eqref{equ-obs} failed, we could choose $\norm{f_n}_2=1$ with
$\int_0^T\norm{\ind_E U_m(t)f_n}_2^2\,dt\to0$ as $n\to\infty$. Since $P_0$ has finite
rank, there is a subsequence such that $P_0f_{n_k}$ converges in $L^2(\R)$.
Applying \eqref{equ-tr} to $f_{n_k}-f_{n_\ell}$ gives
\[
\begin{aligned}
\alpha_T\norm{f_{n_k}-f_{n_\ell}}_2^2
&\le 2N_a\int_0^T\bigl(\norm{\ind_E U_m(t)f_{n_k}}_2^2+
                           \norm{\ind_E U_m(t)f_{n_\ell}}_2^2\bigr)\,dt\\
&\quad+\tfrac T2\norm{P_0f_{n_k}-P_0f_{n_\ell}}_2^2\longrightarrow0
\qquad\text{as }k,\ell\to\infty.
\end{aligned}
\]
Thus $f_{n_k}$ is Cauchy in $L^2(\R)$. Let $f$ be its limit in this space. Then
$\norm f_2=\lim_{k\to\infty}\norm{f_{n_k}}_2=1$ and
\[
\int_0^T\norm{\ind_E U_m(t)f}_2^2\,dt
\le 2T\norm{f-f_{n_k}}_2^2+2\int_0^T\norm{\ind_E U_m(t)f_{n_k}}_2^2\,dt
\longrightarrow0\quad\text{as }k\to\infty.
\]
By strong continuity, $\ind_E U_m(t)f=0$ for every $0<t<T$.

Set $f=\sum_n c_n\varphi_n$. For $y\in L^2(\R)$ supported in $E$ up to a null set,
let $a_n=c_n\int_E\varphi_n(x)\overline{y(x)}\,dx$. By Cauchy--Schwarz and Bessel's inequality,
\[
\sum_n|a_n|\le
\Big(\sum_n|c_n|^2\Big)^{1/2}
\Big(\sum_n\Big|\int_\R\varphi_n\overline y\,dx\Big|^2\Big)^{1/2}
\le\norm f_2\norm y_2.
\]
Define $F_y(z)=\sum_na_ne^{-iz\lambda_n}$ for $\Imm z\le0$.
Since $|e^{-iz\lambda_n}|=e^{\lambda_n\Imm z}\le1$ there, the series converges
uniformly on the closed lower half-plane and its sum is holomorphic in the open lower half-plane.
For $0<t<T$,
\[
F_y(t)=\int_E(U_m(t)f)(x)\overline{y(x)}\,dx=0.
\]
By continuity up to $(0,T)$, Schwarz reflection extends $F_y$ holomorphically
across this interval. The identity theorem then gives $F_y(z)=0$ for $\Imm z<0$.
In particular, for each $s>0$ take $y=\ind_E e^{-sH_m}f$. Then
\[
0=F_y(-is)=\int_E |(e^{-sH_m}f)(x)|^2\,dx
=\norm{\ind_E e^{-sH_m}f}_2^2.
\]
Since $\norm f_2=1$, there is a smallest $j$ with $c_j\ne0$.
Simplicity of the spectrum and orthogonality give
\[
\norm{e^{s\lambda_j}e^{-sH_m}f-c_j\varphi_j}_2^2
=\sum_{n>j}|c_n|^2e^{-2s(\lambda_n-\lambda_j)}
\le e^{-2s(\lambda_{j+1}-\lambda_j)}\norm f_2^2\longrightarrow0
\qquad\text{as }s\to\infty.
\]
Hence $c_j\ind_E\varphi_j=\lim_{s\to\infty}e^{s\lambda_j}\ind_E e^{-sH_m}f=0$
in $L^2(\R)$. Since $c_j\ne0$ and $d_m(E,a)>0$ implies $|E|>0$, this contradicts
the eigenfunction zero-set property in Section~\ref{sec2-pre}.

\end{proof}

\section{Necessary conditions and completion of the proof}
\subsection{Gaussian packets and necessary density conditions}
For $1/2<m\le1$, we retain the quadratic term in the expansion of
$V(x)=|x|^{2m}$ about the classical trajectory to construct a moving Gaussian packet.
For the propagation of coherent states by quadratic approximation, see
\cite{Hagedorn,CR}. For related Gaussian-wavepacket constructions and
observability estimates for time-dependent quadratic Hamiltonians, under
restrictions on the initial data, see Waters~\cite{Waters}.

\begin{lemma}\label{lem-g}
Let $c\in C([0,T];[0,1])$.
The initial value problem
\[
\begin{cases}
i\partial_tz=-\partial_y^2z+c(t)y^2z, &(t,y)\in[0,T]\times\R,\\
z(0,y)=g(y):=\pi^{-1/4}e^{-y^2/2}, &y\in\R,
\end{cases}
\]
has a normalized Gaussian solution satisfying
\[
 |z(t,y)|^2\le C_Te^{-b_Ty^2},\qquad
 C_T=\frac{3e^{4T}}{2\sqrt\pi},\qquad b_T=\frac1{9e^{8T}}>0.
\]
Its spatial Schwartz seminorms are bounded uniformly over such coefficients $c$.
\end{lemma}
\begin{proof}
Consider the initial value problem
\[
v''(t)+4c(t)v(t)=0,\qquad v(0)=1,\qquad v'(0)=2i.
\]
Since $c$ is continuous, it has a unique solution $v\in C^2([0,T];\mathbb C)$.
Since $c$ is real,
\[
\frac d{dt}\Imm(v'\overline v)
=\Imm(-4c|v|^2+|v'|^2)=0.
\]
Hence $\Imm(v'\overline v)=2$ and $v(t)\ne0$ on $[0,T]$.
Integrating the initial value problem gives
\[
v(t)=1+\int_0^t v'(\tau)\,d\tau,\qquad
v'(t)=2i-4\int_0^t c(\tau)v(\tau)\,d\tau.
\]
For $M(t)=|v(t)|+|v'(t)|$, the bound $0\le c\le1$ gives
\[
M(t)\le3+\int_0^t (|v'(\tau)|+4|v(\tau)|)\,d\tau
\le3+4\int_0^t M(\tau)\,d\tau.
\]
By Gronwall's inequality, $|v(t)|+|v'(t)|\le3e^{4t}$. Since
$2=\Imm(v'(t)\overline{v(t)})\le |v'(t)|\,|v(t)|$, we have $2/(3e^{4T})\le |v(t)|\le 3e^{4T}$.
Put
\[
\alpha(t)=\frac{v'(t)}{2iv(t)},\qquad
z(t,y)=\pi^{-1/4}v(t)^{-1/2}e^{-\alpha(t)y^2/2},
\]
where the continuous branch is chosen so that $v(0)^{-1/2}=1$.
Set $A(t)=\pi^{-1/4}v(t)^{-1/2}$. Then
\[
 \Ree\alpha=|v|^{-2},\qquad \alpha'=2i(c-\alpha^2),\qquad iA'/A=\alpha .
\]
Substitution gives $i\partial_t z=-\partial_y^2 z+c(t)y^2z$,
and $\alpha(0)=1$ gives $z(0,y)=g(y)$.
With $r=y/|v(t)|$ we obtain
\[
\int_\R |z(t,y)|^2\,dy
=\frac{1}{\sqrt\pi\,|v(t)|}\int_\R e^{-y^2/|v(t)|^2}\,dy
=\frac1{\sqrt\pi}\int_\R e^{-r^2}\,dr=1.
\]
In particular, for $0\le t\le T$,
\[
|z(t,y)|^2\le\frac{3e^{4T}}{2\sqrt\pi}
\exp\!\left(-\frac{y^2}{9e^{8T}}\right),\quad
|\alpha(t)|=\frac{|v'(t)|}{2|v(t)|}\le\frac94e^{8T},\quad
\Ree\alpha(t)\ge\frac19e^{-8T}.
\]
Boundedness of $\alpha$ and a positive lower bound for $\Ree\alpha$
give uniform spatial Schwartz seminorms, since each derivative is a
polynomial times the same Gaussian.
\end{proof}

\begin{proposition}\label{prop-nec}
Let $L_m=2m+1$.
For $1/2<m<1$, observation in any time $T>0$ implies
\begin{equation}\label{equ-ns}
d_m(E,L_mT^2)>0.
\end{equation}
For $m=1$, the same implication holds for $0<T\leq1/4$.
\end{proposition}
\begin{proof}
\Needspace{4\baselineskip}\medskip\noindent\textbf{Step 1. Approximation by a moving Gaussian packet.}\par\nopagebreak[4]
Start the classical trajectory at rest at $R>0$:
\begin{equation}\label{equ-traj}
q_R'=2p_R,\qquad p_R'=-V'(q_R),\qquad q_R(0)=R,\quad p_R(0)=0.
\end{equation}
The vector field is smooth for $q_R>0$, giving a unique local solution.
Energy conservation gives $p_R^2+V(q_R)=R^{2m}$, so the solution continues
as long as $R/2\le q_R\le R$. We claim that, for all sufficiently large $R$ and every $t\in[0,T]$,
\begin{equation}\label{equ-tb}
4m2^{-s}tR^s\leq-q_R'(t)\leq4mtR^s,\qquad
0\leq R-q_R(t)\leq2mt^2R^s.
\end{equation}
While $R/2\le q_R(t)\le R$, \eqref{equ-traj} gives
\[
4m2^{-s}R^s\leq -q_R''(t)\leq4mR^s.
\]
Integrating from $0$ to $t$, and using $q_R'(0)=0$, yields
\[
4m2^{-s}tR^s\leq-q_R'(t)\leq4mtR^s.
\]
Integrating again gives
\[
2m2^{-s}t^2R^s\leq R-q_R(t)\leq2mt^2R^s.
\]
The trajectory moves left and cannot exit upwards through $R$.
If it first reached $R/2$ at $\tau\leq T$, then
\[
\frac R2=R-q_R(\tau)\leq2mT^2R^s<R/2.
\]
For $s<1$ the last strict inequality holds for all sufficiently large
$R$; for $m=1$ it follows from $T\leq1/4$.
Thus $R/2<q_R(t)\le R$ on $[0,T]$, and \eqref{equ-tb} follows.
For large $R$ we have $q_R\ge1$. Set
\[
c_R(t)=\tfrac12V''(q_R(t))=m(2m-1)q_R(t)^{2m-2}\in[0,1].
\]
Let $z_R(t,y)$ be the Gaussian solution of Lemma~\ref{lem-g} with coefficient $c_R$.
Then $|z_R(t,y)|^2\le C_Te^{-b_Ty^2}$ for all $0\le t\le T$ and $y\in\R$.
Define
\[
S_R(t)=\int_0^t(p_R(\tau)^2-V(q_R(\tau)))\,d\tau,\qquad
\Theta_R(t,x)=p_R(t)(x-q_R(t))+S_R(t),
\]
and $u_R^{\rm app}(t,x)=e^{i\Theta_R(t,x)}z_R(t,x-q_R(t))$.
For brevity, suppress the subscript $R$ and set $y=x-q(t)$. Then
$\Theta_t=p'y-pq'+S'$ and $\Theta_x=p$, so
\begin{align*}
e^{-i\Theta}(i\partial_t-H_m)u^{\rm app}
&=iz_t+z_{yy}+i(2p-q')z_y\\
&\quad+[-p'y+pq'-S'-p^2-V(q+y)]z.
\end{align*}
Substituting $q'=2p$, $p'=-V'(q)$, $S'=p^2-V(q)$ and
$iz_t+z_{yy}=\tfrac12V''(q)y^2z$ gives
\begin{equation}\label{equ-prs}
\begin{split}
(i\partial_t-H_m)u_R^{\rm app}
&=\mathcal R_R(t,x-q_R)u_R^{\rm app},\\
\mathcal R_R(t,y)
&=V(q_R)+V'(q_R)y+\tfrac12V''(q_R)y^2-V(q_R+y).
\end{split}
\end{equation}
If $|y|\le R/4$, then $q_R(t)+\theta y\ge R/2-|y|\ge R/4$ for $0\le\theta\le1$.
Taylor's formula gives
\[
\mathcal R_R(t,y)
=-\frac{y^3}{2}\int_0^1(1-\theta)^2
V^{(3)}(q_R(t)+\theta y)\,d\theta.
\]
There $|V^{(3)}|\leq C_mR^{2m-3}$, and $|z_R(t,y)|^2\le C_Te^{-b_Ty^2}$ gives
\[
\norm{\ind_{\{|y|\leq R/4\}}\mathcal R_Rz_R}_2
\leq C_mR^{2m-3}
\left(\int_\R|y|^6|z_R(t,y)|^2\,dy\right)^{1/2}
\leq C_{m,T}R^{2m-3}.
\]
On $|y|>R/4$, use $|q_R+y|^{2m}\le 2^{2m-1}(|q_R|^{2m}+|y|^{2m})$,
$R\ge1$ and $2m\le2$ to obtain
\[
|\mathcal R_R(t,y)|
\leq C_m(R^{2m}+R^s|y|+R^{2m-2}y^2+|y|^{2m})
\leq C_m(R^2+y^2).
\]
Consequently,
\[
\norm{\ind_{\{|y|>R/4\}}\mathcal R_Rz_R}_2^2
\leq C_{m,T}\int_{|y|>R/4}(R^4+y^4)e^{-b_Ty^2}\,dy
\leq C_{m,T}R^4e^{-b_TR^2/32}.
\]
Thus for $m<1$,
\begin{equation}\label{equ-rs}
\sup_{0\leq t\leq T}\norm{\mathcal R_R(t,\cdot)z_R(t,\cdot)}_2
\leq C_{m,T}(R^{2m-3}+R^2e^{-b_TR^2/64})\longrightarrow0
\quad\text{as }R\to\infty.
\end{equation}
At $m=1$, $V(x)=x^2$, so $\mathcal R_R(t,y)=0$ for every $t,y$.
Put $F_R(t,x)=\mathcal R_R(t,x-q_R(t))u_R^{\rm app}(t,x)$ and
$w_R(t)=U_m(t)g(\cdot-R)-u_R^{\rm app}(t)$. The two functions have the same initial value $g(x-R)$ at $t=0$, and $(i\partial_t-H_m)w_R=-F_R$. Duhamel's formula gives
\[
w_R(t)=i\int_0^t U_m(t-\tau)F_R(\tau)\,d\tau.
\]
By unitarity, followed by the change of variables $y=x-q_R(\tau)$ and
$|e^{i\Theta_R}|=1$,
\[
\begin{aligned}
\norm{w_R(t)}_2
&=\left\|\int_0^t U_m(t-\tau)F_R(\tau)\,d\tau\right\|_2
\le\int_0^t\norm{U_m(t-\tau)F_R(\tau)}_2\,d\tau\\
&=\int_0^t\norm{F_R(\tau)}_2\,d\tau
=\int_0^t\norm{\mathcal R_R(\tau,\cdot)z_R(\tau,\cdot)}_2\,d\tau.
\end{aligned}
\]
For $m<1$, \eqref{equ-rs} gives
\begin{equation}\label{equ-dp}
\eta_R:=\sup_{0\le t\le T}\norm{w_R(t)}_2
\le T\sup_{0\le t\le T}\norm{\mathcal R_R(t,\cdot)z_R(t,\cdot)}_2\longrightarrow0
\quad\text{as }R\to\infty.
\end{equation}
At $m=1$, $F_R=0$, so $U_1(t)g(\cdot-R)=u_R^{\rm app}(t)$ and $\eta_R=0$.

\Needspace{4\baselineskip}\medskip\noindent\textbf{Step 2. A lower density bound from observability.}\par\nopagebreak[4]
Set $C=C_{\obs,m}(E,T)<\infty$.
For $u_R(t)=U_m(t)g(\cdot-R)$, observability and normalization give
\[
1=\norm{g(\cdot-R)}_2^2\le C\int_0^T\!\int_E |u_R(t,x)|^2\,dx\,dt .
\]
By \eqref{equ-dp},
\[
\int_0^T\!\int_E |u_R-u_R^{\rm app}|^2\,dx\,dt
\le\int_0^T\norm{u_R(t)-u_R^{\rm app}(t)}_2^2\,dt\le T\eta_R^2 .
\]
For sufficiently large $R$, $\sqrt T\,\eta_R\le(2\sqrt{C})^{-1}$. Hence
\[
\norm{u_R^{\rm app}}_{L^2((0,T)\times E)}
\ge \norm{u_R}_{L^2((0,T)\times E)}-\sqrt T\,\eta_R
\ge C^{-1/2}-\sqrt T\,\eta_R\ge(2\sqrt{C})^{-1}.
\]
Set $d_T=(4C)^{-1}$, $\delta=\min\{T/2,d_T/4\}$ and
$I_R=I_{m,L_mT^2}(R)=[R-L_mT^2R^s,R+L_mT^2R^s]$.
Since $|u_R^{\rm app}(t,x)|=|z_R(t,x-q_R(t))|$,
\begin{equation}\label{equ-mm}
\begin{aligned}
d_T &\le \int_0^T\!\int_E |z_R(t,x-q_R(t))|^2\,dx\,dt\\
&=\underbrace{\int_0^\delta\!\int_E |z_R(t,x-q_R(t))|^2\,dx\,dt}_{\Pi_1}
 +\underbrace{\int_\delta^T\!\int_{E\setminus I_R}|z_R(t,x-q_R(t))|^2\,dx\,dt}_{\Pi_2}\\
&\quad +\underbrace{\int_\delta^T\!\int_{E\cap I_R}|z_R(t,x-q_R(t))|^2\,dx\,dt}_{\Pi_3}.
\end{aligned}
\end{equation}

\Needspace{4\baselineskip}\medskip\noindent\textbf{\boldmath Estimate of $\Pi_1$.}\par\nopagebreak[4]
Since $\norm{z_R(t,\cdot)}_2=1$, we have $\Pi_1\le\delta\le d_T/4$.

\Needspace{4\baselineskip}\medskip\noindent\textbf{\boldmath Estimate of $\Pi_2$.}\par\nopagebreak[4]
For $x\notin I_R$, \eqref{equ-tb} gives
\[
 |x-q_R(t)|\ge |x-R|-|R-q_R(t)|\ge (L_m-2m)T^2R^s=T^2R^s.
\]
With $y=x-q_R(t)$, Lemma~\ref{lem-g} therefore gives
\[
 \Pi_2\le (T-\delta) C_T\int_{|y|\ge T^2R^s}e^{-b_Ty^2}\,dy
 \longrightarrow0\qquad\text{as }R\to\infty,
\]
since $s>0$. In particular, $\Pi_2\le d_T/4$ for all sufficiently large $R$.

\Needspace{4\baselineskip}\medskip\noindent\textbf{\boldmath Estimate of $\Pi_3$.}\par\nopagebreak[4]
By \eqref{equ-mm}, $\Pi_3\ge d_T-\Pi_1-\Pi_2\ge d_T/2$.
On $[\delta,T]$, \eqref{equ-tb} gives
$-q_R'(t)\ge 4m2^{-s}\delta R^s$. For fixed $x$, the change of variables
$y=x-q_R(t)$ gives
\[
\begin{aligned}
\int_\delta^T|z_R(t,x-q_R(t))|^2\,dt
&\le \int_{x-q_R(\delta)}^{x-q_R(T)}
 \frac{C_Te^{-b_Ty^2}}{-q_R'(t(y))}\,dy\\
&\le \frac{C_T\sqrt{\pi/b_T}}{4m2^{-s}\delta R^s}.
\end{aligned}
\]
Integrating over $x\in E\cap I_R$ yields
\[
 \frac{d_T}{2}\le\Pi_3\le
 \frac{C_T\sqrt{\pi/b_T}}{4m2^{-s}\delta R^s}|E\cap I_R|.
\]
Since $|I_R|=2L_mT^2R^s$,
\begin{equation}\label{equ-nd}
\frac{|E\cap I_{m,L_mT^2}(R)|}{2L_mT^2R^s}
\geq\frac{d_T(4m2^{-s})\delta}{4L_mT^2 C_T\sqrt{\pi/b_T}}>0
\qquad(R\text{ sufficiently large}).
\end{equation}
By reflection, testing with $g(\cdot+R)$ gives the same bound at $-R$.
This proves \eqref{equ-ns}.
\end{proof}

\subsection{Proof of Theorem~\ref{thm-d}}
\begin{proof}
For $1/2<m<1$, if $d_m(E,a)>0$ for some $a>0$, Proposition~\ref{prop-suf}
gives observability for $T>\pi\sqrt{A_ma}$. Conversely, observability in
time $T$ implies $d_m(E,(2m+1)T^2)>0$ by Proposition~\ref{prop-nec}.
At $m=1$, positive density at some scale is equivalent to weak thickness by
Lemma~\ref{lem-wt}, and weak thickness characterizes observability in some positive
time by Huang--Wang--Wang~\cite[Theorem~1.2(i)--(ii)]{HWW}.
This proves \eqref{equ-st}.

For \eqref{equ-et}, let $1/2<m\le1$ and suppose $d_m(E,a)>0$ for every $a>0$.
Given $T>0$, choose $0<a<T^2/(\pi^2A_m)$, with $a\le a_0$ if $m=1$.
Proposition~\ref{prop-suf} gives observability in time $T$.
Conversely, assume that $E$ is observable for every $T>0$. For any $a>0$
choose $0<T<\min\{1/4,\sqrt{a/(2m+1)}\}$ and set $b=(2m+1)T^2<a$.
Proposition~\ref{prop-nec} gives $d_m(E,b)>0$, and Lemma~\ref{lem-ds} yields
\[
 d_m(E,a)\ge\frac{b}{a}\,d_m(E,b)>0 .
\]

Finally, let $1/2<m<1$. By Proposition~\ref{prop-nec}, every admissible
time $T$ satisfies $\kappa_m(E)\le(2m+1)T^2$. Taking the infimum gives
$\Tmin(E)\ge\sqrt{\kappa_m(E)/(2m+1)}$.
If $\kappa_m(E)<\infty$, then for every $a>\kappa_m(E)$, Lemma~\ref{lem-ds}
gives $d_m(E,a)>0$, and Proposition~\ref{prop-suf} gives
$\Tmin(E)\le\pi\sqrt{A_ma}$. Letting $a\downarrow\kappa_m(E)$ proves the upper
bound in \eqref{equ-tc}. If $\kappa_m(E)=\infty$, no finite time
is admissible and $\Tmin(E)=\infty$.
\end{proof}

\phantomsection

\Needspace{4\baselineskip}
\medskip
\noindent{\fontsize{9}{11}\selectfont\scshape School of Mathematics and Statistics, Central South University, Changsha, Hunan 410083, P.R. China}\\
\textit{Email address:} \href{mailto:yanchuanwang@csu.edu.cn}{\texttt{yanchuanwang@csu.edu.cn}}

\end{document}